\documentclass[
reqno,
11pt,
oneside,
a4paper
]{article}

\usepackage[ngerman, english]{babel}
\usepackage{enumerate}
\usepackage[normalem]{ulem}
\usepackage[babel,french=guillemets,german=swiss]{csquotes}
\usepackage[center,font=small]{caption}
\usepackage{cite}
\usepackage{bbm}
\usepackage[raggedright]{titlesec}
\usepackage{xcolor}

\titleformat{\section}{\normalfont\large\bfseries}{\thesection}{1em}{}
\titleformat{\subsection}{\normalfont\bfseries}{\thesubsection}{1em}{}

\usepackage{tocbasic}
\DeclareTOCStyleEntry[
beforeskip=.2em plus 1pt,
]{tocline}{section}

\usepackage[%
left = 2.7cm,
right = 2.7cm,
top = 2cm,
bottom = 2cm,
footskip = 1cm,
]{geometry}

\usepackage{color}
\definecolor{LinkColor}{rgb}{0,0,1}
\definecolor{LinkColor2}{rgb}{0,0.5,0}
\definecolor{lbcolor}{rgb}{0.85,0.85,0.85}
\definecolor{FrameColor}{rgb}{0.85,0.85,0.85}
\definecolor{rosso}{rgb}{0.8,0,0}
\definecolor{lightgray}{rgb}{0.5,0.5,0.5}
\definecolor{violet}{rgb}{0.65,0,0.65}
\definecolor{darkgreen}{rgb}{0,0.5,0}

\usepackage{enumitem}

\usepackage{graphicx}
\usepackage{placeins}
\usepackage{overpic}

\usepackage{graphicx,tikz}
\usetikzlibrary{patterns,calc,hobby,positioning} 
\usepackage[]{empheq} 

\usepackage{amsmath}
\usepackage{amssymb}
\usepackage{dsfont}
\usepackage{mathrsfs}
\usepackage{amsthm}

\usepackage[%
pdftitle={Titel},%
pdfauthor={Autor},%
pdfcreator={LaTeX, LaTeX with hyperref and KOMA-Script},
pdfsubject={Betreff}, 
pdfkeywords={Keywords}
]{hyperref} 

\hypersetup{%
	colorlinks	=true,
	linkcolor	=LinkColor,%
	anchorcolor	=LinkColor,%
	citecolor	=LinkColor2,%
	filecolor	=LinkColor,%
	menucolor	=LinkColor,%
	urlcolor	=LinkColor,%
}
\usepackage{marginnote}

\usepackage{verbatim}

\newtheorem{theorem}{Theorem}[section]

\theoremstyle{definition}
\newtheorem{remark}[theorem]{Remark}

\allowdisplaybreaks[4]



\makeatletter
\renewcommand\paragraph{\@startsection{paragraph}{4}{\z@}%
	{1ex \@plus1ex \@minus.2ex}%
	{-1em}%
	{\normalfont\normalsize\bfseries}}
\renewcommand\subparagraph{\@startsection{paragraph}{4}{\z@}%
	{1ex \@plus1ex \@minus.2ex}%
	{-1em}%
	{\normalfont\normalsize\itshape}}
\makeatother

\newcommand{\cA}{\mathcal{A}}

\newcommand{\cD}{\mathcal{D}}
\newcommand{\cE}{\mathcal{E}}

\newcommand{\cV}{\mathcal{V}}

\newcommand{\n}{{\mathbf n}}
\newcommand{\uu}{{\mathbf u}}
\newcommand{\vv}{{\mathbf v}}

\newcommand{\x}{{\mathbf x}}
\newcommand{\y}{{\mathbf y}}

\newcommand{\D}{\mathbf{D}}
\newcommand{\I}{\mathbf{I}}

\newcommand{\J}{\mathbf{J}}
\newcommand{\K}{\mathbf{K}}

\newcommand{\T}{\mathbf{T}}

\newcommand{\bS}{\mathbf{S}}

\newcommand{\R}{\mathbb R}
\newcommand{\N}{\mathbb N}

\newcommand{\abs}[1]{\left| #1 \right|}

\newcommand{\mom}{m_\Omega}
\newcommand{\mga}{m_\Gamma}

\newcommand{\trho}{\tilde{\rho}}
\newcommand{\tsigma}{\tilde{\sigma}}

\newcommand{\tJ}{\tilde{\mathbf J}}
\newcommand{\tK}{\tilde{\mathbf K}}

\newcommand{\tT}{\tilde{\mathbf T}}

\newcommand{\Om}{\Omega(t)}
\newcommand{\Ga}{\Gamma(t)}

\newcommand{\ST}{S_T}
\newcommand{\QT}{Q_T}

\newcommand{\intO}{\int_{\Om}}

\newcommand{\intG}{\int_{\Ga}}

\newcommand{\eps}{\varepsilon}

\newcommand{\dx}{\;\mathrm d\x}
\newcommand{\dX}{\;\mathrm dX}

\newcommand{\dt}{\;\mathrm dt}

\newcommand{\dH}{\;\mathrm d\mathcal{H}}
\newcommand{\dHH}[1]{\dH^{ #1 }}

\newcommand{\h}{\mathds{h}}

\newcommand{\ddt}{\frac{\mathrm d}{\mathrm dt}}

\newcommand{\del}{\partial}
\newcommand{\delt}{\partial_{t}}
\newcommand{\deltb}{\partial_{t}^{\bullet}}
\newcommand{\deltc}{\partial_{t}^{\circ}}
\newcommand{\deltn}{\partial_{t}^{\scalebox{0.5}{$\square$}}}

\newcommand{\deln}{\partial_\n}

\newcommand{\delph}{\partial_{\phi}}
\newcommand{\delgph}{\partial_{\Grad\phi}}
\newcommand{\delps}{\partial_{\psi}}
\newcommand{\delgps}{\partial_{\Gradg\psi}}

\newcommand{\Grad}{\nabla}
\newcommand{\Lap}{\Delta}
\newcommand{\Div}{\textnormal{div}\,}
\newcommand{\Gradg}{\nabla_\Gamma}
\newcommand{\Lapg}{\Delta_\Gamma}
\newcommand{\Divg}{\textnormal{div}_\Gamma}

\newcommand{\ov}{\overline}
\newcommand{\suchthat}{\;\ifnum\currentgrouptype=16 \middle\fi|\;}

\newcommand{\tin}{\quad \text{in }}
\newcommand{\ton}{\quad \text{on }}

\begin{document}
	
	%
	%


    \title{\bfseries\Large
		A Thermodynamically Consistent Free Boundary Model \\
        for Two-Phase Flows in an Evolving Domain \\
		with Bulk-Surface Interaction
	}
	
	\author{
		Patrik Knopf \footnotemark[1]
		\and Yadong Liu \footnotemark[2]}
	
	\date{ }
	
	\maketitle
	
	\renewcommand{\thefootnote}{\fnsymbol{footnote}}
	
	\footnotetext[1]{
		Fakult\"at f\"ur Mathematik, 
		Universit\"at Regensburg, 
		93053 Regensburg, 
		Germany
		\tt(%
		\href{mailto:patrik.knopf@ur.de}{patrik.knopf@ur.de}%
		).
	}
	\footnotetext[2]{
		School of Mathematical Sciences, Ministry of Education Key Laboratory of NSLSCS, and Key Laboratory of Jiangsu Provincial Universities of FDMTA, Nanjing Normal University, Nanjing 210023, P.~R.~China 
		\tt(%
		\href{mailto:ydliu@njnu.edu.cn}{ydliu@njnu.edu.cn}%
		).
	}
	
    
	%
	%
	
	\begin{small}
		\begin{center}
			\textbf{Abstract}
		\end{center}
		We derive a thermodynamically consistent model, which describes the time evolution of a two-phase flow in an evolving domain. The movement of the free boundary of the domain is driven by the velocity field of the mixture in the bulk, which is determined by a Navier--Stokes equation. In order to take interactions between bulk and boundary into account, we further consider two materials on the boundary, which may be the same or different materials as those in the bulk. The bulk and the surface materials are represented by respective phase-fields, whose time evolution is described by a bulk-surface convective Cahn--Hilliard equation. This approach allows for a transfer of material between bulk and surface as well as variable contact angles between the diffuse interface in the bulk and the boundary of the domain. To provide a more accurate description of the corresponding contact line motion, we include a generalized Navier slip boundary condition on the velocity field. Based on local mass balance laws, we derive our model from scratch in two different ways: by the Lagrange Multiplier Approach and (in the case of matched densities and no mass flux between bulk and surface) by the Energetic Variational Approach. We further show that our model generalizes previous models from the literature, which can be recovered from our system by either dropping the dynamic boundary conditions or assuming a static boundary of the domain.
		\\[1.5ex]
		\\[1.5ex]	
	\end{small}

	
	\begin{footnotesize}
		\setcounter{tocdepth}{2}
		\hypersetup{linkcolor=black}
		\tableofcontents
	\end{footnotesize}

	\setlength\parskip{0.5ex}
	\allowdisplaybreaks
	\numberwithin{equation}{section}
	\renewcommand{\thefootnote}{\arabic{footnote}}
	
	\newpage
	
	\section{Introduction}\label{sec:introduction}

    The description of two-phase flows is an important topic in modern continuum fluid dynamics, which has widespread applications in biology, chemistry, and engineering.
    To represent the interface separating the different components, two fundamental approaches have been developed, namely \textit{sharp-interface methods} and \textit{diffuse interface methods}. We refer to \cite{Du2020,Abels2018,Giga2018,PS2016} for a comparison of the two approaches.

    In sharp-interface models, the interface is represented by an evolving hypersurface contained in the surrounding domain. The time evolution of the interface is then described by a free boundary problem. 

    In diffuse interface models, the interface between two materials is represented by a thin layer whose thickness is proportional to a small parameter $\eps>0$. The location of the two components is usually represented by an order parameter, the so-called \textit{phase-field}, which describes the difference of the local concentrations (or volume fractions) of the materials.
    This means that the phase-field attains values close to $-1$ or $1$ in the regions where the single fluids are present. At the diffuse interface between fluids, the phase-field is expected to exhibit a continuous transition between the values representing the pure phases. In normal direction across the diffuse interface, a certain transition profile can usually be observed. The main advantage of this approach is that the time evolution of the phase-field can be described via an Eulerian formulation as a PDE on the whole domain. This avoids directly tracking the interface between the phases. If a proper scaling with respect to the interface parameter $\eps$ is used, the corresponding sharp-interface model can be recovered (in most cases at least formally) as the {\it sharp-interface limit} (see, e.g., \cite{Abels2018}).  

    One of the most fundamental models for describing the motion of two viscous and incompressible fluids with equal densities is the \textit{Model H}.
    It was first proposed in \cite{Hohenberg1977} and later rigorously derived in \cite{Gurtin1996}. The system consists of an incompressible Navier--Stokes equation that is coupled to a convective Cahn--Hilliard equation.
    The fundamental assumption in the derivation of the Model H is that the individual constant densities of fluids are equal. They are therefore referred to as \textit{matched densities}. This directly entails that the density of the whole mixture is constant.  
    
    In the past two decades, various works have proposed generalizations of the Model H  
    describing incompressible mixtures with \textit{unmatched} densities and/or compressible two-phase flows. 
    We refer to
    \cite{Abels2012,Lowengrub1998,Boyer2002,Ding2007, Shen2013,Shokrpour2018,Giga2018,Heida2012,Freistuhler2017}
    to at least mention some of them.
    
    One of the most popular models for viscous incompressible two-phase flows with unmatched densities is the one derived by Abels, Garcke and Grün in \cite{Abels2012}. In view of the initials of the authors, it is often referred to as the \textit{AGG Model}. It was shown in \cite{Abels2012} that by including an additional flux term in the Navier--Stokes equation, the model can be made thermodynamically consistent even if the constant individual densities of fluids are \textit{unmatched}. Analytical results for the AGG model can be found, e.g., in \cite{Abels2013,Abels2013a,Abels2021,Abels2024,Abels2024a,Giorgini2021,Giorgini2022}. An alternative derivation of the AGG model via mixture theory was recently developed in \cite{tenEikelder2023}.

    However, as the classical AGG model is considered subject to a no-slip boundary condition for the velocity field as well as homogeneous Neumann boundary conditions for the Cahn--Hilliard quantities (i.e., the phase-field and the chemical potential), it still has some drawbacks if a precise description of the behavior of fluids close to the boundary is required. 
    In the Cahn--Hilliard equation, the homogeneous Neumann condition on the phase-field enforces the diffuse interface to always intersect the boundary at a ninety degree contact angle, which is very unrealistic in many applications. The homogeneous Neumann boundary condition on the chemical potential entails the conservation of mass within the considered domain, which makes a description of transfer of material between bulk an boundary (e.g., absorption processes) impossible. 
    These issues have also been noticed by physicists 
    (see, e.g., \cite{Binder1991,Fischer1997,Fischer1998,Kenzler2001}),
    who proposed the inclusion of an additional surface free energy to describe short-range interactions between bulk and boundary more precisely. 
    
    Due to these observations, several types of dynamic boundary conditions for the Cahn--Hilliard equation have been proposed and analyzed in the literature. 
    Especially in recent times, dynamic boundary conditions that also exhibit a Cahn--Hilliard type structure have become very popular because they allow for a quantification of the mass transfer between bulk and boundary. For brevity, we refer to such models as \textit{bulk-surface Cahn--Hilliard} models.
    We refer to  \cite{Garcke2020,Garcke2022,Colli2020,Colli2022,Colli2022a,Fukao2021,Miranville2020} for analytical results and to \cite{Metzger2021,Metzger2023,Harder2022,Meng2023,Bao2021,Bao2021a,Bullerjahn2024,Bullerjahn2025} for a numerical investigation of the bulk-surface Cahn--Hilliard system. A nonlocal bulk-surface Cahn--Hilliard system was analyzed in \cite{Knopf2021b,Lv2025}. 
	A convective variant of the bulk-surface Cahn--Hilliard system with a prescribed velocity field was analyzed, for example, in \cite{Colli2018,Colli2019,Gilardi2019,Giorgini2025,Knopf2024,Knopf2024a}.
    The role of dynamic boundary conditions in Cahn--Hilliard type models will be discussed in more detail in Subsection~\ref{SUBSEC:DYN}.
    For more information on the Cahn--Hilliard equation with classical homogeneous Neumann boundary conditions or with dynamic boundary conditions, we refer to the recent review paper \cite{Wu2022} as well as the book \cite{Miranville-Book}.    

    Moreover, as shown for example in \cite{Qian2006,Seppecher1996}, a no-slip boundary condition on the velocity field is not very well suited for the description of moving contact line phenomena. This is because convective effects occurring at the boundary (or close to the boundary) are neglected if the velocity field is assumed to be identically zero at the surface. Therefore, if a no-slip boundary condition is imposed, the motion of the contact line is primarily driven by diffusion or convection processes away from the boundary. To overcome this issue, a generalized Navier slip boundary condition, which provides a better description of moving contact lines, was derived in \cite{Qian2006}. 

    Based on these considerations, Navier--Stokes--Cahn--Hilliard models with dynamic boundary conditions on the Cahn--Hilliard subsystem as well as a generalized Navier slip boundary condition for the Navier--Stokes equation have been proposed and investigated in the literature. 
    We refer to \cite{Gal2016,Gal2019}, where a second-order (Allen--Cahn type) dynamic boundary condition was used for the Cahn--Hilliard subsystem. Here, the Allen--Cahn type dynamic boundary condition can be interpreted as a parabolic relaxation of a transport-type dynamic boundary condition related to contact angle dynamics (see, e.g., \cite{Yue2020}).
    Instead, in \cite{Giorgini2023}, a fourth-order (Cahn--Hilliard type) dynamic boundary condition is considered, which allows for a better description of the transfer of material between bulk and boundary.
    This model was further analyzed in \cite{Gal2024}. In these models, the Cahn--Hilliard quantities at the boundary enter the generalized Navier slip boundary condition via a forcing term. We also refer to \cite{Feireisl2010,Cherfils2019}, where \textit{compressible} two-phase flows were considered. In \cite{Feireisl2010} a compressible Navier--Stokes--Allen--Cahn model was analyzed and in \cite{Cherfils2019}, a compressible Navier--Stokes--Cahn--Hilliard model with dynamic boundary conditions was investigated.
    Moreover, models for two-phase flows with \textit{phase-transition} were studied in \cite{Aki2014,Abels2025a}.

\begin{figure}[ht]
    	\centering
    	\includegraphics[width=0.6\textwidth]{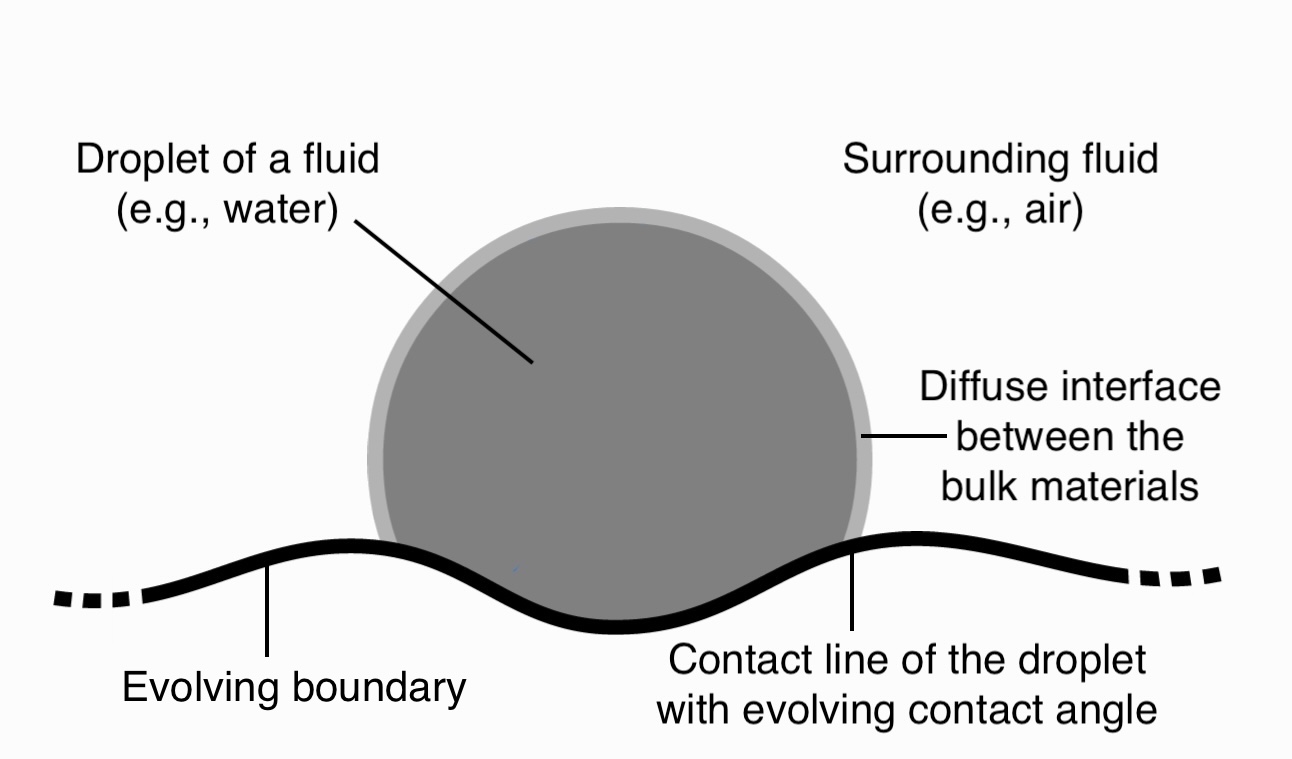}
    	\caption{Droplet on a moving surface.\\ 
        (Only the relevant part of the evolving domain $\Omega(\cdot)$ and its boundary is displayed.)}
        \label{FIG:DROP}
    \end{figure}

    We point out that all the models discussed so far are each considered in a fixed bounded domain that does not change during the time evolution. In some applications, however, the domain in which a two-phase flow is considered will evolve over the course of time. A simple example is a droplet (surrounded by a further substance such as air) that is placed on a deformable substrate, for example, a deformable sheet (see Figure~\ref{FIG:DROP}). 
    
\begin{figure}[ht]
        \vspace{3ex}
    	\centering
    	\includegraphics[width=0.75\textwidth]{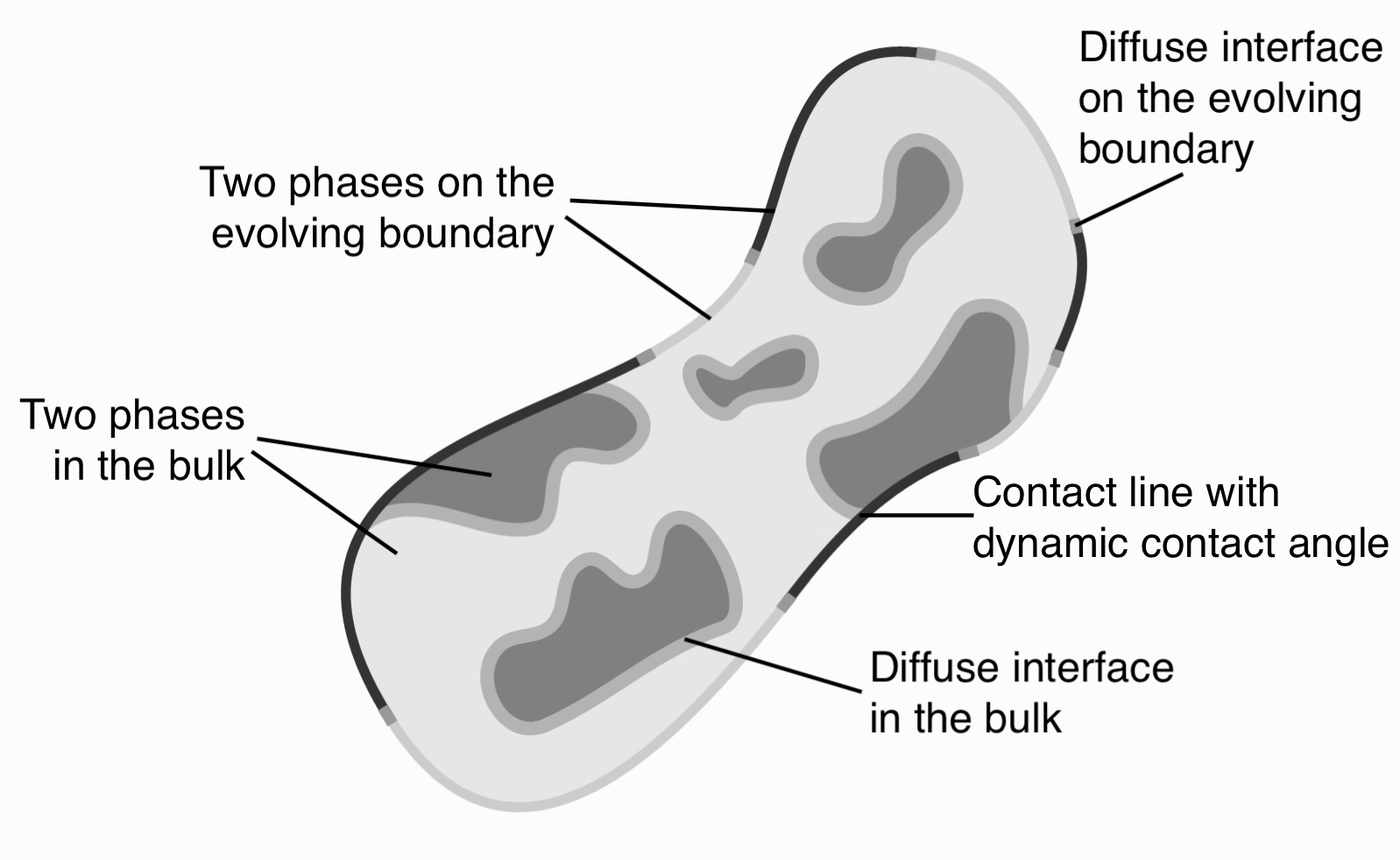}
    	\caption{Evolving domain (e.g., describing a bacterium) with\\ 
        two phases in the bulk and two phases on the evolving boundary.}
        \label{FIG:BACT}
    \end{figure}

    A more complex example is the motion of a single biological cell, for instance, a bacterium (cf.~Figure~\ref{FIG:BACT}). In this case, the bacterium is represented by an evolving domain and its cell wall, represented by the boundary of the domain, is an evolving surface. Different substances inside the bacteria (e.g., cell organelles) can be described as a two-phase (or even multi-phase) flow. Also the materials on the cell wall (such as proteins or polysaccharides) can be interpreted as a mixture of different substances whose time evolution is of interest. 

    In this paper, for the aforementioned reasons, we derive a new model, which generalizes the two-phase flow model derived in \cite{Giorgini2023} to the framework of an evolving domain with a free boundary. 
    Here, a key assumption is that the evolution of our domain is driven by the velocity field of the mixture contained in the domain. This means that the normal velocity of the free boundary is equal to the normal component of the velocity field at the boundary. 

    \subsection{Formulation of our model} \label{SUBSEC:FORMMOD}
    
	To formulate our system of equations, we consider an evolving domain $\{\Omega(t)\}_{t\in[0,T]}$ in $\R^d$ (with $d\in \N$, $T>0$). We set $\Gamma(t)\coloneqq\partial\Omega(t)$ for all $t\in[0,T]$. This means that $\{\Gamma(t)\}_{t\in[0,T]}$ is a closed, oriented evolving hypersurface in $\R^d$. For brevity, we use the notation
	\begin{equation*}
		\QT \coloneqq \bigcup_{0 < t <T} \{t\} \times \Om 
		\quad\text{and}\quad 
		S_T \coloneqq \bigcup_{0 < t <T} \{t\} \times \Ga 
	\end{equation*}
	We further assume that the evolving domain $\QT$ is transported by a sufficiently regular vector field $\uu:\overline{Q_T}\to\R^d$, which is supposed to be the \textit{volume-averaged velocity field} of the two materials that are contained inside the domain.
	This entails that $S_T$ is transported by the trace vector field $\uu\vert_{S_T}:S_T\to\R^d$. Since $S_T$ is orientable, there exists an evolving unit normal vector field $\n$ and the normal velocity $\cV$ of $S_T$ is thus given by
	\begin{align*}
		\uu \cdot \n = \cV \quad\text{on $S_T$}.
	\end{align*}
    For a mathematically precise formulation of the notion of evolving domains and surfaces that are transported by a velocity field, we refer to \cite{BDGP-book,BGN-book}.
    
	In order to describe the time evolution of the velocity field $\uu$, we derive the following thermodynamically consistent model:%
    \begin{subequations}
		\label{eqs:NSCH}
        \begin{alignat}{3}
			\label{eqs:NSCH:1}
			& \delt \big(\rho(\phi) \uu\big) 
			+ \Div \big( \uu \otimes (\rho(\phi)\uu + \J) \big) 
			= \Div(\T),
			\quad \Div \uu = 0
			&& \tin \QT, 
			\\
			\label{eqs:NSCH:2}
			& \deltb \phi 
			= \Div (m_\Omega(\phi) \Grad \mu)
			&& \tin \QT,
			\\
			\label{eqs:NSCH:3}
			&\mu 
			= - \varepsilon \Lap \phi + \varepsilon^{-1} F'(\phi)  && \tin \QT,
			\\
			\label{eqs:NSCH:4}
			& \deltc \psi  
			= \Divg (m_\Gamma(\psi) \Gradg \theta) - \beta m_\Omega(\phi) \deln \mu  && \ton S_T,
			\\
			\label{eqs:NSCH:5}
			& \theta 
			= - \delta \Lapg \psi + \delta^{-1} G'(\psi)
			+ \alpha \varepsilon \deln \phi  && \ton S_T, 
			\\
			\label{eqs:NSCH:6}
			&K \varepsilon \deln \phi
			= \alpha \psi - \phi, 
			\quad 
			L m_\Omega(\phi) \deln \mu
			= \beta \theta - \mu, 
			\quad K,L \in [0,\infty], 
			&& \ton S_T, 
			\\
			\label{eqs:NSCH:7}
			& \big[\T\n  + \Gradg q  + \gamma_\tau  \uu \big]_\tau
			= \Big[ 
                \frac{1}{2} \, (\J\cdot\n) \, \uu
                - \delta\, \Divg(\Gradg\psi \otimes \Gradg\psi)
            \Big]_\tau
			&& \ton S_T, 
			\\
			\label{eqs:NSCH:8}
			& \T\n\cdot\n + \gamma_\n \cV 
			= \frac{1}{2} \, (\J\cdot\n) \cV 
                - q H 
			&& \ton S_T, \\
			\label{eqs:NSCH:9}
			&  \Divg \uu = 0,  
            \quad \uu \cdot \n = \cV
			&& \ton S_T.
        \end{alignat}
	\end{subequations}
	This system can be understood as a generalized version of the model introduced in \cite{Giorgini2023}, adapted to the situation of an evolving domain with free boundary. 
    An equivalent reformulation of this model will be discussed in Subsection~\ref{SUBSEC:ALTFOR}.

    In the above equations and throughout this paper, $\Grad$ denotes the gradient, $\Gradg$ denotes the tangential gradient on the surface, $\Div$ denotes the divergence, and $\Divg$ denotes the surface divergence. For (sufficently regular) functions $\mathbf{w}:\QT\to \R^d$ and $\mathbf{M}:\QT\to \R^{d\times d}$, we use the conventions
    \begin{equation*}
        [\Grad\mathbf{w}]_{ij} = \partial_j \mathbf{w}_i
        \quad\text{and}\quad
        [\Div \mathbf{M}]_i = \sum_{k=1}^d \partial_k \mathbf{M}_{ik}
    \end{equation*}
    for all $i,j\in\{1,...,d\}$. This means that for any vector-valued function, $\Grad$ provides its Jacobian matrix.
    For functions on the surface $S_T$, an analogous convention is employed.
    Moreover, for any (column) vectors $\mathbf{v}, \mathbf{w}\in\R^d$, their tensor product is given by $\mathbf{v} \otimes \mathbf{w} = \mathbf{v}\mathbf{w}^T$.
	
	In system \eqref{eqs:NSCH}, the volume-averaged velocity field of the two fluids is represented by the vector-valued function $\uu:\QT\to\R^d$. The scalar functions $\phi:\QT\to\R$, $\mu:\QT\to\R$, $\psi:S_T\to\R$ and $\theta:S_T\to\R$ denote the \textit{bulk phase-field}, the \textit{bulk chemical potential}, the \textit{surface phase-field}, and the  \textit{surface chemical potential}, respectively.
    The parameters $\eps>0$ and $\delta>0$ are related to the thickness of the diffuse interface in the bulk and on the evolving boundary, respectively.
	It is a key feature of the model \eqref{eqs:NSCH} that the evolving domain $\QT$ is transported by the velocity field $\uu$. This is reflected in the boundary condition $\eqref{eqs:NSCH:9}_2$ and, therefore, system \eqref{eqs:NSCH} is a \textit{free boundary problem}.
	
	The time evolution of $\uu$ is described by the \textit{incompressible Navier--Stokes equations} \eqref{eqs:NSCH:1}. Here, $\rho(\phi)$ stands for the phase-field dependent \textit{density}, $\J$ represents the \textit{mass flux} related to interfacial motion, and $\T$ is the \textit{stress tensor}. They are given by the formulae
	\begin{align}
		\label{eq:rho}
		\rho(\phi) &\coloneqq \frac{1}{2} \trho_2 (1+\phi) + \frac{1}{2} \trho_1 (1-\phi),
		\\
		\label{eq:J}
		\J &\coloneqq - \frac{1}{2} (\trho_2 - \trho_1) \mom(\phi) \Grad\mu,
		\\
		\label{eq:T}
		\T &\coloneqq  \bS - p\I - \eps \Grad\phi \otimes \Grad\phi
		\quad\text{with}\quad
		\bS \coloneqq 2\nu(\phi) \, \D\uu,
	\end{align}
	where the constants $\trho_1,\trho_2>0$ are the individual densities of fluids, $\nu(\phi)$ is the phase-field dependent \textit{viscosity}, and
	$$\D\uu \coloneqq \frac{1}{2} \big( (\Grad \uu) + (\Grad \uu)^T \big)$$
	denotes the \textit{symmetric gradient} of $\uu$. 
    The scalar function $p: Q_T \to \R$ denotes the \textit{bulk pressure}. It can be interpreted as a Lagrange multiplier to account for the incompressibility condition $\eqref{eqs:NSCH:1}_2$. The matrix $\bS$ is referred to as the \textit{viscous stress tensor}.
	
	The evolution of the phase-fields $\phi$ and $\psi$ is described by the \textit{bulk-surface convective Cahn--Hilliard system} \eqref{eqs:NSCH:2}--\eqref{eqs:NSCH:6}. 
	Here, the \textit{bulk Cahn--Hilliard equation} \eqref{eqs:NSCH:2}--\eqref{eqs:NSCH:3} is coupled to the \textit{surface Cahn--Hilliard equation} \eqref{eqs:NSCH:4}--\eqref{eqs:NSCH:5}, which can be interpreted as a dynamic boundary condition. To state the Cahn--Hilliard equations, we are using the following notation. For scalar functions $f:\QT\to\R$ and $g:S_T\to\R$, we write
	\begin{equation}
		\label{DEF:MATDER}
		\deltb f \coloneqq \delt f + \uu\cdot \nabla f,
		\qquad 
		\deltc g  
        \coloneqq \delt \tilde g + \uu \cdot \Grad \tilde g
        = \delt \tilde g + \uu_\tau \cdot \Gradg \tilde g + \cV \n\cdot \Grad \tilde g
	\end{equation}
	to denote the \textit{material derivatives} of $f$ and $g$, respectively. Here, $\tilde g$ denotes an arbitrary (sufficiently regular) extension of $g$ onto a neighborhood of $S_T$. In fact, the expression $\deltc g$ does actually not depend on the choice of this extension (see, e.g., \cite[Remark~2.7.2]{BDGP-book}).
    The phase-field dependent functions $m_\Omega(\phi)$ and $m_\Gamma(\psi)$ represent the \textit{mobilities} of the bulk- and the surface materials, respectively. Moreover, $F'$ and $G'$ are the derivatives of \textit{double-well potentials} $F$ and $G$, respectively.
	The bulk and surface quantities of the Cahn--Hilliard equations are further coupled through the boundary conditions in \eqref{eqs:NSCH:6}, which are to be understood as follows:%
	\begin{subequations}
		\label{CH:CPL}
		\begin{alignat}{2}
			\label{CH:CPL:1}
			&\begin{cases}
				\phi = \alpha \psi  &\text{if $K=0$},\\
				K \varepsilon \deln \phi = \alpha \psi - \phi &\text{if $K\in (0,\infty)$},\\
				\deln \phi = 0 &\text{if $K=\infty$}
			\end{cases}
			&&\quad \ton S_T, 
			\\[1ex]
			\label{CH:CPL:2}
			&\begin{cases}
				\mu = \beta \theta  &\text{if $L=0$},\\
				L m_\Omega(\phi) \deln \mu = \beta \theta - \mu &\text{if $L\in (0,\infty)$},\\
				m_\Omega(\phi) \deln \mu = 0 &\text{if $L=\infty$}
			\end{cases}
			&&\quad \ton S_T.
		\end{alignat}
	\end{subequations}
	Here, $\alpha,\beta\in\R$ are prescribed parameters. 
	The parameters $K,L\in [0,\infty]$ are used to distinguish different cases, each resulting in a certain solution behaviour related to a corresponding physical phenomenon. They will be explained in more detail in Subsection~\ref{SUBSEC:DYN}. 
    If $K=0$, we additionally assume that $\alpha\neq 0$.

    Furthermore, on the surface $S_T$, the tangential component of the velocity field is governed by the \textit{generalized Navier slip boundary condition} \eqref{eqs:NSCH:7} with a slip parameter $\gamma_\tau$ that may depend on the materials involved (i.e., the phase-fields).
	Here, for any vector field $\mathbf{w}$, the notation $\mathbf{w}_\tau$ indicates that only the tangential component of $\mathbf{w}$ is considered.
    Condition $\eqref{eqs:NSCH:9}_1$ is imposed as we assume the evolving boundary to be inextensible.
    In particular, for simplicity, surface elasticity is neglected in our model. This means that
    $\Gamma(t)$ is to be interpreted as a massless interface without any stored surface elastic energy.
    The function $q: S_T \to \R$ appears as a Lagrange multiplier accounting for the inextensibility condition $\eqref{eqs:NSCH:9}_1$.  
	The normal component of $\uu$ on the boundary is determined by the relation \eqref{eqs:NSCH:8}, which involves the \textit{mean curvature} $H = - \frac{1}{d}\Divg\n$.
    In particular, the term $qH$ represents the impact of the geometry of the boundary to the surface normal traction.
	Moreover, $\gamma_\n>0$ is a prescribed scalar function that may depend on the phase-fields. It is related to energy dissipation due to boundary deformation (cf.~\eqref{diss}). 

    The energy functional associated with system \eqref{eqs:NSCH} is given by
    \begin{equation}
        \label{DEF:E}
        \begin{aligned}
            \cE(\uu,\phi,\psi) & \coloneqq 
            \intO \frac{\rho(\phi)}{2} \abs{\uu}^2 \dx
            + \intO \Big(
            \frac{\varepsilon}{2} \abs{\Grad \phi}^2 
            + \frac{1}{\varepsilon} F(\phi)
            \Big) \dx \\
            & \qquad + \intG \Big(
            \frac{\delta}{2} \abs{\Gradg \psi}^2 
            + \frac{1}{\delta} G(\psi)
            + \frac{1}{2}\mathds{h}(K) \abs{\alpha \psi - \phi}^2
            \Big) \dHH{d-1},
        \end{aligned}
    \end{equation}
    where the function 
    \begin{equation}
        \label{DEF:H}
        \mathds{h}:[0,\infty] \to \R, \quad
        \mathds{h}(r) \coloneqq
        \begin{cases}
            r^{-1} &\text{if $r\in(0,\infty)$},\\
            0 &\text{if $r\in\{0,\infty\}$}
        \end{cases}
    \end{equation}
    is used to account for the different choices of the parameter $K$. In \eqref{DEF:E}, the first summand on the right-hand side represents the \textit{kinetic energy}. The second summand is the \textit{bulk free energy} and the third summand is the \textit{surface free energy}. They are both of \textit{Ginzburg--Landau type} and the free energy densities crucially depend on the potentials $F$ and $G$, respectively, which are assumed to be double-well shaped. This means that they are coercive with two local minimum points close to $-1$ and $1$ (as these values represent the pure phases of the materials) and a local maximum point in between (usually at zero).
    A typical choice for $F$ and $G$ is the so-called \textit{Flory--Huggins potential}
    \begin{align}\label{DEF:W:LOG}
        W_{\textup{log}}(s) \coloneqq \frac{\Theta}{2}\big[(1+s)\ln(1+s) + (1-s)\ln(1-s)\big] -\frac{\Theta_c}{2}s^2, \quad s\in(-1,1)
    \end{align}
    for $0 < \Theta < \Theta_c$, which can be derived from Boltzmann kinetics. It is sometimes simply referred to as the \textit{logarithmic potential}. 
    Since $W_{\textup{log}}'(s) \to \pm\infty$ as $s\to\pm 1$, the potential $W_{\textup{log}}$ is a so-called \textit{singular potential}. This property ensures that for any solution of system~\eqref{eqs:NSCH} with $F=G=W_{\textup{log}}$, the phase-fields $\phi$ and $\psi$ attain their values within the physically relevant range $[-1,1]$. An important consequence is that due to $\abs{\phi}\le 1$, we have
    \begin{equation*}
        0 < \min\{\trho_1,\trho_2\} \leq \rho(\phi) \leq \max\{\trho_1,\trho_2\} < \infty, 
    \end{equation*}
    which means that the density $\rho(\phi)$ remains uniformly positive and finite.
    Furthermore, our model has the following physically relevant properties:
	\begin{itemize}[leftmargin=*, topsep=0em, partopsep=0em, parsep=1ex, itemsep=0ex]
		\item \textbf{Energy dissipation.}
		Sufficiently regular solutions of system \eqref{eqs:NSCH} satisfy the \textit{energy-dissipation law}
		\begin{equation}
			\label{diss}
			\begin{aligned}
				& \ddt \cE(\uu,\phi,\psi)
				+ \intO 2 \nu(\phi) \abs{\D \uu}^2 \dx
				+ \intG \gamma_{\tau} \abs{\uu_{\tau}}^2 \dHH{d-1} 
				\\&\quad 
				+ \intG \gamma_{\n} \abs{\uu \cdot \n}^2 \dHH{d-1} 
				+ \intO m_{\Omega}(\phi) \abs{\Grad \mu}^2 \dx
				\\&\quad 
				+ \intG m_{\Gamma}(\psi) \abs{\Gradg \theta}^2 \dHH{d-1}
				+ \intG \mathds{h}(L) \abs{\beta \theta - \mu}^2 \dHH{d-1}
				= 0
			\end{aligned}
		\end{equation}
		for all $t\in [0,T]$. In this sense, our model can be considered \textit{thermodynamically consistent}.
		\item \textbf{Conservation of mass.}
		Sufficiently regular solutions of system \eqref{eqs:NSCH} fulfill the \textit{mass conservation law} 
		\begin{align}
			\label{mass}
			&\beta \intO \phi(t) \dx + \intG \psi(t) \dHH{d-1}
			= \beta \int_{\Omega(0)} \phi(0) \dx + \int_{\Gamma(0)} \psi(0) \dHH{d-1}
		\end{align}
		for all $t\in [0,T]$. In the case $L=\infty$, it even holds
		\begin{align}
			\label{mass:inf}
			\intO \phi(t) \dx 
			= \int_{\Omega(0)} \phi(0) \dx
			\quad\text{and}\quad
			\intG \psi(t) \dHH{d-1} 
			= \int_{\Gamma(0)} \psi(0) \dHH{d-1}
		\end{align}
		for all $t\in [0,T]$, so the bulk mass and the surface mass are conserved separately.
		\item \textbf{Conservation of bulk volume and surface area.}
		For any sufficiently regular solution of system \eqref{eqs:NSCH}, the volume of the evolving domain $\Omega(\cdot)$ and the area of the evolving boundary $\Gamma(\cdot)$ are preserved over the course of time.
		Namely, as the bulk fluid is assumed to be incompressible (i.e., $\Div\uu = 0$ in $\QT$) and there is no extra source, it follows that
		\begin{align}
			\label{vol}
			\ddt \mathrm{Vol}_{d}\big(\Omega(t)\big)
			= \ddt \intO 1 \dx 
			= \intO \Div(\uu) \dx 
			= 0 
			\quad\text{for all $t\in[0,T]$}
		\end{align}
		by means of the Reynolds transport theorem (see, e.g., \cite[Theorem~2.11.1]{BDGP-book}).
		Similarly, as we also demand $\Divg\uu = 0$ on $\ST$, the transport theorem for evolving hypersurfaces (see, e.g., \cite[Theorem~2.10.1]{BDGP-book}) yields
		\begin{align}
			\label{area}
			\ddt \mathrm{Vol}_{d-1}\big(\Gamma(t)\big)
			= \ddt \intG 1 \dHH{d-1} 
			= \intG \Divg(\uu) \dHH{d-1} 
			= 0 
			\quad\text{for all $t\in[0,T]$}
		\end{align}
        because of $\eqref{eqs:NSCH:9}_1$.
	\end{itemize}

    \subsection{Further related literature}

    In contrast to our model, where the free boundary is the boundary of the evolving domain, there are also models for two-phase flows, where a free boundary contained in the interior of a static (i.e., time-independent) domain is considered.
    We refer, for example, to \cite{Garcke2014,Barrett2017}. 
    A related sharp interface model with moving contact lines on inextensible elastic sheets was derived in \cite{YZR2023SharpInterface}.
    For the investigation of Cahn--Hilliard type models on evolving surfaces without any bulk contribution, we refer to \cite{Abels2025,Bachini2023,Bachini2023a,Caetano2021,Caetano2023,Elliott2024} as well as the references therein.

    \subsection{The outline of this paper}
    
	The plan of this paper is as follows.
    In Section~\ref{SECT:VAR}, we present an equivalent reformulation of system \eqref{eqs:NSCH}. Moreover, we compare two simplified variants of our model, one without dynamic boundary conditions and one in a domain with static boundary, with previous models proposed in the literature. We further discuss the benefits of dynamic boundary conditions in our model as well as in Cahn--Hilliard type models in general.
	Section~\ref{SECT:MODEL} is then devoted to the derivation of the Navier--Stokes--Cahn--Hilliard model~\eqref{eqs:NSCH}. Starting from the local mass balance laws, we derive the model by two different approaches: the Lagrange Multiplier Approach and, in the case of matched densities and no mass flux between bulk and surface, the Energetic Variational Approach.

    \section{Important features and variants of our model}\label{SECT:VAR} 
	
	\subsection{An alternative reformulation}\label{SUBSEC:ALTFOR} 
    Instead of $p$ and $q$, we can alternatively use the pressures
    \begin{alignat*}{2}
        \ov p &:= p + f
        &&\quad\text{with}\quad
        f \coloneqq 
        \Big( \frac \eps2 \abs{\Grad\phi}^2 +\frac 1\eps F(\phi) \Big),
        \\
        \ov q &:= q + g
        &&\quad\text{with}\quad
        g \coloneqq 
        \Big( \frac \delta 2 \abs{\Gradg\psi}^2 
            + \frac 1\delta G(\psi) 
            + \frac 12 \h(K) (\alpha\psi-\phi)^2 \Big),
    \end{alignat*}
    which involve the free energy densities $f$ and $g$.
    It is straightforward to check that
    \begin{align*}
        &- \Grad p - \eps\, \Div\big(\Grad\phi \otimes \Grad\phi \big)
        \\
        &\quad= \mu \Grad\phi 
            - \Grad p 
            - \Grad \Big( \frac \eps2 \abs{\Grad\phi}^2 +\frac 1\eps F(\phi) \Big) 
        = \mu \Grad\phi 
            - \Grad \ov p
    \end{align*}
    in $\QT$, and
    \begin{align*}
        &\eps\,(\Grad\phi \otimes \Grad\phi)\n  
            - \Gradg q 
            - \delta\, \Divg\big(\Gradg\psi \otimes \Gradg\psi \big)
        \\
        &\quad= - \Gradg q 
            - \Gradg \Big( \frac \delta2 \abs{\Gradg\psi}^2 
                +\frac 1\delta G(\phi)
                + \frac 12\h(K) (\alpha\psi-\phi)^2
            \Big)
        \\
        &\quad\qquad + \theta \Gradg \psi 
            - \eps \deln\phi (\alpha\Gradg \psi - \Grad\phi)   
            + \h(K) (\alpha\psi-\phi) (\alpha\Gradg \psi - \Gradg\phi) 
        \\[1ex]
        &\quad= \theta \Gradg\psi 
            - \Gradg \ov q
            - \eps \deln\phi (\alpha\Gradg \psi - \Grad\phi)   
            + \h(K) (\alpha\psi-\phi) (\alpha\Gradg \psi - \Gradg\phi)
    \end{align*}
    on $\ST$.
    By the definition of $\h(K)$ and boundary condition \eqref{eqs:NSCH:6}, this implies that
    \begin{align*}
        &\Big[\eps\,(\Grad\phi \otimes \Grad\phi)\n  
            - \Gradg q 
            - \delta\, \Divg\big(\Gradg\psi \otimes \Gradg\psi \big) \Big]_\tau
        \\
        &\quad= \theta \Gradg\psi 
            - \Gradg \ov q
            - \big[ \eps \deln\phi - \h(K) (\alpha\psi-\phi) \big]
                (\alpha\Gradg \psi - \Gradg\phi)  
        \\
        &\quad= \theta \Gradg\psi 
            - \Gradg \ov q.
    \end{align*}
    on $\ST$.
    Recalling the definitions of the stress tensor $\T$ and the viscous stress tensor $\bS$ (see~\eqref{eq:T}), system~\eqref{eqs:NSCH} can thus be reformulated as\pagebreak[3]  
    \begin{subequations}
		\label{eqs:NSCH:ALT}
        \begin{alignat}{3}
			\label{eqs:NSCH:ALT:1}
			& \delt \big(\rho(\phi) \uu\big) 
			+ \Div \big( \uu \otimes (\rho(\phi)\uu + \J) \big) 
            - \Div(\bS) + \Grad \ov p
			= \mu \Grad\phi,
			\quad \Div \uu = 0
			&& \tin \QT, 
			\\
			\label{eqs:NSCH:ALT:2}
			& \deltb \phi 
			= \Div (m_\Omega(\phi) \Grad \mu),
			&& \tin \QT,
			\\
			\label{eqs:NSCH:ALT:3}
			&\mu 
			= - \varepsilon \Lap \phi + \varepsilon^{-1} F'(\phi)  && \tin \QT,
			\\
			\label{eqs:NSCH:ALT:4}
			& \deltc \psi  
			= \Divg (m_\Gamma(\psi) \Gradg \theta) - \beta m_\Omega(\phi) \deln \mu  && \ton S_T,
			\\
			\label{eqs:NSCH:ALT:5}
			& \theta 
			= - \delta \Lapg \psi + \delta^{-1} G'(\psi)
			+ \alpha \varepsilon \deln \phi  && \ton S_T, 
			\\
			\label{eqs:NSCH:ALT:6}
			&K \varepsilon \deln \phi
			= \alpha \psi - \phi, 
			\quad 
			L m_\Omega(\phi) \deln \mu
			= \beta \theta - \mu, 
			\quad K,L \in [0,\infty], 
			&& \ton S_T, 
			\\
			\label{eqs:NSCH:ALT:7}
			& \big[\bS\n  + \Gradg \ov q  + \gamma_\tau  \uu \big]_\tau
			= \Big[ 
                \frac{1}{2} \, (\J\cdot\n) \, \uu
                + \theta \Gradg\psi
            \Big]_\tau
			&& \ton S_T, 
			\\
			\label{eqs:NSCH:ALT:8}
			& \bS\n\cdot\n 
                - \ov p
                - \eps (\deln\phi)^2
                + \gamma_\n \cV 
			= \frac{1}{2} \, (\J\cdot\n) \cV 
                - \ov q H + gH - f
			&& \ton S_T, \\
			\label{eqs:NSCH:ALT:9}
			&  \Divg \uu = 0,  
            \quad \uu \cdot \n = \cV
			&& \ton S_T.
        \end{alignat}
	\end{subequations}
    This formulation might be more convenient for future mathematical analysis of this model.
	
	\subsection{Variants of our model without dynamic boundary conditions}\label{SUBSEC:FBAGG}
	
	The standard choice of boundary conditions for the classical Cahn--Hilliard equation, where no surface quantities are considered, are homogeneous Neumann boundary conditions on both $\phi$ and $\mu$. This situation is also included in our model as a special case. 
	By assuming $\psi\equiv -1$ or $\psi\equiv 1$ and choosing $K=L=\infty$, we obtain $\theta\equiv 0$ and therefore, system \eqref{eqs:NSCH} reduces to%
	\begin{subequations}
		\label{eqs:AGGQ}
        \begin{alignat}{3}
			\label{eqs:AGGQ:1}
			& \delt \big(\rho(\phi) \uu\big) 
			+ \Div \big( \uu \otimes (\rho(\phi)\uu + \J) \big) 
			= \Div(\T),
			\quad \Div \uu = 0
			&& \tin \QT, 
			\\
			\label{eqs:AGGQ:2}
			& \deltb \phi
			= \Div (m_\Omega(\phi) \Grad \mu)
			&& \tin \QT,
			\\
			\label{eqs:AGGQ:3}
			&\mu 
			= - \varepsilon \Lap \phi + \varepsilon^{-1} F'(\phi)  && \tin \QT,
			\\
			\label{eqs:AGGQ:6}
			&\deln \phi = 0, 
			\quad 
			\deln \mu = 0 
			&& \ton S_T, 
			\\
			\label{eqs:AGGQ:7}
			& \big[\bS\n  + \Gradg q  + \gamma_\tau  \uu \big]_\tau
			= \Big[ \frac{1}{2} \, (\J\cdot\n) \, \uu\Big]_\tau
			&& \ton S_T, 
			\\
			\label{eqs:AGGQ:8}
			& \T\n\cdot\n + \gamma_\n  \cV 
			= \frac{1}{2} \, (\J\cdot\n) \, \cV 
            - qH
			&& \ton S_T, \\
			\label{eqs:AGGQ:9}
			&  \Divg\uu = 0 , \quad \uu \cdot \n = \cV
			&& \ton S_T,
        \end{alignat}
	\end{subequations} 
	where $\rho(\phi)$, $\J$ and $\T$ are given by \eqref{eq:rho}--\eqref{eq:T}.
    
	As a further simplification, it is possible to neglect the condition $\Divg\uu = 0$ on $\ST$, which also allows us also to omit the surface pressure $q$ as it acts as a Lagrange multiplier for this condition (cf.~Subsection~\ref{SUBSEC:LAG}). 
    This simplification was also employed in \cite{Giorgini2023}, where a variant of our model with a static boundary was derived. There, neither the condition $\Divg\uu = 0$ nor any surface pressure were included in the model. In this way, the system \eqref{eqs:AGGQ} further simplifies to
	\begin{subequations}
		\label{eqs:AGG}
        \begin{alignat}{3}
			\label{eqs:AGG:1}
			& \delt \big(\rho(\phi) \uu\big) 
			+ \Div \big( \uu \otimes (\rho(\phi)\uu + \J) \big) 
			= \Div(\T),
			\quad \Div \uu = 0
			&& \tin \QT, 
			\\
			\label{eqs:AGG:2}
			& \deltb \phi
			= \Div (m_\Omega(\phi) \Grad \mu)
			&& \tin \QT,
			\\
			\label{eqs:AGG:3}
			&\mu 
			= - \varepsilon \Lap \phi + \varepsilon^{-1} F'(\phi)  && \tin \QT,
			\\
			\label{eqs:AGG:6}
			&\deln \phi = 0, 
			\quad 
			\deln \mu = 0
			&& \ton S_T, 
			\\
			\label{eqs:AGG:7}
			& \big[\bS\n + \gamma_\tau  \uu \big]_\tau
			= \Big[ \frac{1}{2} \, (\J\cdot\n) \, \uu\Big]_\tau
			&& \ton S_T, 
			\\
			\label{eqs:AGG:8}
			& \T\n\cdot\n + \gamma_\n  \cV 
			= \frac{1}{2} \, (\J\cdot\n) \, \cV 
			&& \ton S_T, 
            \\
			\label{eqs:AGG:9}
			&\uu \cdot \n = \cV
			&& \ton S_T.
        \end{alignat}
	\end{subequations} 
	This version can be regarded as a variant of the Abels--Garcke--Grün model (cf.~\cite{Abels2012}) in an evolving domain with free boundary.

    \subsection{The role of dynamic boundary conditions in our model}\label{SUBSEC:DYN}
    
	As already discussed in the literature (see, e.g., \cite{Giorgini2023}), the homogeneous Neumann boundary conditions in \eqref{eqs:AGG:6} are quite restrictive. On the one hand, the condition $\deln \phi = 0$ on $S_T$ means that the contact angle between the diffuse interface separating fluids and the boundary $\Gamma(\cdot)$ is always exactly ninety degrees. Of course, this is unrealistic in many applications, where the contact angle is expected to change dynamically over the course of time. On the other hand, the condition $\deln \mu = 0$ on $S_T$ implies that no transfer of material between bulk and surface can occur (cf.~\eqref{mass:inf}). Therefore, absorption or adsorption processes cannot be described. 
	
	These issues are overcome in model \eqref{eqs:NSCH} by means of the bulk-surface convective Cahn--Hilliard subsystem \eqref{eqs:NSCH:2}--\eqref{eqs:NSCH:6}. 
	At least for the non-convective bulk-surface Cahn--Hilliard system (i.e., $\uu\equiv 0$), it has been observed in several numerical simulations that a transfer of material between bulk and surface occurs if $L\in [0,\infty)$ and that the contact angle between the diffuse interface and the boundary changes dynamically over the course of time (see, e.g., \cite{Knopf2021a,Bao2021a,Metzger2023,Bullerjahn2024,Bullerjahn2025}). 
	Unfortunately, the sharp-interface limit of such bulk-surface Cahn--Hilliard models, which would provide an explicit law for the behavior of the contact angle, is still not very well understood. 
    In \cite{Metzger2025,Zhang2023}, the relation between the Cahn--Hilliard equation with dynamic boundary conditions and a generalized Young's formula for the contact angle was discussed at least on a heuristic level.
	
	In a fixed domain, the bulk-surface Cahn--Hilliard system subject to the coupling conditions \eqref{CH:CPL} with $K,L\in[0,\infty]$ has already been studied extensively in the literature. In particular, the coupling conditions \eqref{CH:CPL} with $K,L\in[0,\infty]$ were introduced and investigated in the following papers:
	\begin{itemize}[leftmargin=*, topsep=0em, partopsep=0em, parsep=1ex, itemsep=0ex]
		\item For $K=L=0$, the non-convective bulk-surface Cahn--Hilliard equation was proposed in \cite{Gal2006,Goldstein2011}. Since $L=0$, the chemical potentials are assumed to always remain in a chemical equilibrium, and a rapid transfer of material between bulk and surface is expected to occur.
		\item For $K=0$ and $L=\infty$, the non-convective bulk-surface Cahn--Hilliard equation was derived in \cite{Liu2019} by the Energetic Variational Approach. In this case, the mass flux between bulk and surface is zero. Consequently, no transfer of material between bulk and surface will occur. 
		\item The case $K\in (0,\infty)$ and $L=\infty$ was first proposed an analyzed in \cite{Knopf2020}. It was shown that solutions of this model converge to solutions of the system with $K=0$ and $L=\infty$ as the parameter $K$ is sent to zero.
		\item For $K=0$ and $L\in (0,\infty)$, the non-convective bulk-surface Cahn--Hilliard equation was proposed in \cite{Knopf2021a} to interpolate between the models associated with $K=0$, $L=0$ and ${K=0}$, ${L=\infty}$. In this case, a transfer of material between bulk and surface will occur, and the number $L^{-1}$ is related to the rate of absorption (or adsorption) of bulk material by the boundary (cf.~\cite{Knopf2021a}).
		\item In \cite{Knopf2024,Knopf2024a,Giorgini2025}, the convective bulk-surface Cahn--Hilliard equation with a given velocity field was studied for all cases $K,L\in[0,\infty]$ in a unified framework. 
		Of course, these results also include the non-convective case by simply choosing $\uu\equiv 0$ for the prescribed velocity field. 
	\end{itemize}
    Moreover, the coupling parameter $\alpha$ in \eqref{CH:CPL} is related to the scaling or mismatch between the phase-fields $\phi$ and $\psi$. In this sense, the most relevant 
    choices are $\alpha=1$ and $\alpha=-1$, but also other choices are possible. 
    The coupling parameter $\beta$ is related to the transfer of material between bulk and boundary.
    More precisely, it represents the fraction of bulk material that is effectively transformed into surface material when reaching the boundary. In this sense, the most relevant choices are $\beta\in [-1,1]$.
	
	\subsection{Variants of our model with a static boundary}
	
	It is also possible to consider a static boundary, which is not allowed to move along with the velocity field. To implement this, we assume $\Omega(\cdot)\equiv \Omega$, $\Gamma(\cdot)\equiv\Gamma$ and $\cV = \uu\cdot\n = 0$ on $\ST$. Then, the simplified model can be expressed as
	\begin{subequations}
		\label{eqs:GKQ}
        \begin{alignat}{3}
			\label{eqs:GKQ:1}
			& \delt \big(\rho(\phi) \uu\big) 
			+ \Div \big( \uu \otimes (\rho(\phi)\uu + \J) \big) 
			= \Div(\T),
			\quad \Div \uu = 0,
			&& \tin \QT, 
			\\
			\label{eqs:GKQ:2}
			& \delt \phi + \Div(\phi\uu)
			= \Div (m_\Omega(\phi) \Grad \mu),
			&& \tin \QT,
			\\
			\label{eqs:GKQ:3}
			&\mu 
			= - \varepsilon \Lap \phi + \varepsilon^{-1} F'(\phi),  && \tin \QT,
			\\
			\label{eqs:GKQ:4}
			& \delt \psi + \Divg(\psi \uu_\tau) 
			= \Divg (m_\Gamma(\psi) \Gradg \theta) - \beta m_\Omega(\phi) \deln \mu,  && \ton S_T,
			\\
			\label{eqs:GKQ:5}
			& \theta 
			= - \delta \Lapg \psi + \delta^{-1} G'(\psi)
			+ \alpha \varepsilon \deln \phi , && \ton S_T, 
			\\
			\label{eqs:GKQ:6}
			&K \varepsilon \deln \phi
			= \alpha \psi - \phi, 
			\quad 
			L m_\Omega(\phi) \deln \mu
			= \beta \theta - \mu, 
			\quad K,L \in [0,\infty], 
			&& \ton S_T, 
			\\
			\label{eqs:GKQ:7}
			& \big[\bS\n  + \Gradg q  + \gamma_\tau  \uu \big]_\tau
			= \Big[ \theta \Gradg \psi + \frac{1}{2} \, (\J\cdot\n) \, \uu\Big]_\tau,
			&& \ton S_T, 
			\\
			\label{eqs:GKQ:9}
			&  \Divg\uu = 0, \quad \uu\cdot\n = 0, 
			&& \ton S_T,
        \end{alignat}
	\end{subequations} 
	where $\QT$ and $\ST$ are to be interpreted as $\QT = \Omega\times (0,T)$ and $\ST = \Gamma\times (0,T)$, respectively. This model can be understood as a variant of the Navier--Stokes--Cahn--Hilliard model derived in \cite{Giorgini2023}, which further includes the condition $\Divg\uu = 0$ on $\ST$ along with the surface pressure $q$ acting as a Lagrange multiplier for this condition. We point out that it becomes clear by the model derivation (see Section~\ref{SECT:MODEL}) that the boundary condition \eqref{eqs:NSCH:8} drops out when the condition $\uu\cdot\n = 0$ on $\ST$ is demanded. 
	
	If we further neglect the relation $\Divg\uu = 0$ on $\ST$ and set $q\equiv 0$, we arrive exactly at the model that was derived in \cite{Giorgini2023}.
	A Cahn--Hilliard--Brinkman variant of this model, in which the Navier--Stokes equation is replaced by a Brinkman equation to describe flows through porous media, was investigated in \cite{Colli2024}.

	\section{Model derivation}\label{SECT:MODEL}
	
	In this section we will derive the Navier--Stokes--Cahn--Hilliard system with dynamic boundary conditions \eqref{eqs:NSCH}.
    Let $\{\Omega(t)\}_{t\in[0,T]}$, $\{\Gamma(t)\}_{t\in[0,T]}$, $\QT$, and $\ST$ be as introduced in Subsection~\ref{SUBSEC:FORMMOD}.
	Recall that the evolving domain $\QT$ is transported by a vector field $\uu$, which is supposed to be the \textit{volume-averaged velocity field} of the two fluids that are considered.
	Therefore, $S_T$ is transported by the vector field $\uu\vert_{S_T}$. 
    Moreover, since $S_T$ is orientable, there exists an evolving unit normal vector field $\n\in C(S_T)$ and the normal velocity $\cV$ of $S_T$ is given by
	\begin{align}
		\label{BC:NORVEL}
		\uu \cdot \n = \cV \quad\text{on $S_T$}.
	\end{align}
    
	In the following, we write $\uu_\tau$ and $\uu_\n$ to denote the tangential component and the normal component of $\uu$, respectively.
	We further assume that each bulk material has a constant individual density $\trho_i$, $i=1,2$
	and each boundary material has a constant individual density $\tsigma_i$, $i=1,2$.
	
	\subsection{Considerations based on local mass balance laws}
	
	The mass densities $\rho_i:\QT\to\R$, $i=1,2$, of the two fluids in the bulk are supposed to satisfy the following mass balance law:
	\begin{align}
		\label{BL:MASS:B}
		\delt \rho_i + \Div \, \widehat\J_i = 0,
		\quad i=1,2,
		\quad\text{in $\QT$.}
	\end{align}
	Here, $\widehat\J_i$ are the mass fluxes corresponding to the motion of the two materials in the bulk. 
	In our model, we also want to account for a transfer of material between bulk and surface. This is, for instance, needed in order to describe absorption processes. 
	Therefore, we also consider functions $\sigma_i:\ST\to\R$, $i=1,2$, which represent the mass densities of the substrates on the free boundary. In general, the materials on the surface might differ from those in the bulk. This is the case, for example, if the materials are transformed by chemical reactions occurring at the boundary. Therefore, we interpret the densities $\sigma_i$ as independent functions that are not necessarily related to the traces of the densities $\rho_i$ (up to a constant carrying physical units), respectively.
	We assume that the following mass balance law holds on the free boundary:
	\begin{align}
		\label{BL:MASS:S}
		\deltn \sigma_i + \Divg \, \widehat\K_i = \widehat\J_{\Gamma,i} \cdot \n,
		\quad i=1,2,
		\quad\text{on $\ST$.}
	\end{align}
	In this relation, $\widehat\K_i$ are (not necessarily tangential) vector fields representing the mass fluxes on the boundary. 
	Moreover, $\widehat\J_{\Gamma,i}$ are (not necessarily tangential) vector fields describing the transfer of material between bulk and boundary. 
	Here and in what follows, the expression 
	\begin{align}
		\label{DEF:NTD}
		\deltn h 
		= \deltc h - \uu_\tau \cdot \Gradg h
        = \delt \tilde h + \cV\n \cdot \Grad \tilde h
	\end{align}
	denotes the \textit{normal time derivative} of any (sufficiently regular) function $h:\ST\to\R$. 
    Here, $\tilde h$ denotes an arbitrary (sufficiently regular) extension of $h$ onto a neighborhood of $S_T$.
	In \eqref{BL:MASS:S}, the normal time derivative needs to be used instead of $\delt\sigma_i$ to account for the normal velocity of the evolving hypersurface $\ST$. 
	In order to ensure the local preservation of surface area (cf.~\eqref{PDE:DIVV:S} below), we additionally postulate that
	\begin{align}
		\label{CONS:VOL}
		\frac{\widehat\J_{\Gamma,1}}{\tsigma_1}\cdot \n + \frac{\widehat\J_{\Gamma,2}}{\tsigma_2}\cdot \n = 0
		\quad\text{on $\ST$}.
	\end{align}	
	We next assume that the motion of our materials is described by the individual velocity fields $\uu_i:\QT\cup\ST\to\R^d$, $i=1,2$.  
	In view of the mass fluxes and the mass densities, we have the following relations:
	\begin{align*}
		\uu_i = \frac{\widehat\J_i}{\rho_i} \quad\text{in $\QT$},
		\qquad
		\uu_i = \frac{\widehat \K_i}{\sigma_i} \quad\text{on $\ST$},
		\qquad
		i=1,2.
	\end{align*}
	Consequently, the mass balances \eqref{BL:MASS:B} and \eqref{BL:MASS:S} can be rewritten as
	\begin{alignat}{2}
		\label{BL:MASS:B'}
		\delt \rho_i + \Div(\rho_i\uu_i) &= 0,
		&&\quad i=1,2,
		\quad\text{in $\QT$,}
		\\
		\label{BL:MASS:S'}
		\deltn \sigma_i  + \Divg(\sigma_i\uu_{i}) &= \widehat\J_{\Gamma,i} \cdot \n,
		&&\quad i=1,2,
		\quad\text{on $\ST$.}
	\end{alignat}
	We now introduce $\phi_i:\QT\to\R$ and $\psi_i:\ST\to\R$, $i=1,2$, to denote the volume fractions of the two fluids in the bulk and the two substrates on the boundary, respectively. 
	Assuming that the bulk and surface materials each have a constant density $\trho_i$, $\tilde{\sigma}_i$, $i=1,2$, the volume fractions can be identified as
	\begin{align}
		\label{DEF:PHIPSI*}
		\phi_i = \frac{\rho_i}{\trho_i} \quad\text{in $\QT$,}
		\qquad
		\psi_i = \frac{\sigma_i}{\tilde{\sigma}_i} \quad\text{on $\ST$,}
		\qquad
		i=1,2.
	\end{align}
	Under the assumption that the excess volume is zero, we have
	\begin{align}
		\label{CONS:PHIPSI}
		\phi_1 + \phi_2 = 1 \quad\text{in $\QT$,}
		\qquad
		\psi_1 + \psi_2 = 1 \quad\text{on $\ST$.}
	\end{align}
	We further define the order parameters
	\begin{align}
		\label{DEF:PHIPSI}
		\phi = \phi_2 - \phi_1 \quad\text{in $\QT$,}
		\qquad
		\psi = \psi_2 - \psi_1 \quad\text{on $\ST$,}
	\end{align}
	which are referred to as the \textit{phase-fields}.
	Using the phase-fields $\phi$ and $\psi$, the velocity field $\uu$ can be decomposed as
	\begin{alignat}{2}
		\label{DEF:VAVG:B}
		\uu &= \phi_1 \uu_1 + \phi_2 \uu_2
		\quad\text{in $\QT$,}
		\\
		\label{DEF:VAVG:S}
		\uu &= \psi_1 \uu_{1} + \psi_2 \uu_{2}
		\quad\text{on $\ST$.}
	\end{alignat}
	We now introduce 
	\begin{alignat}{2}
		\label{Ji}
		\J_i &= \widehat\J_i - \rho_i \uu,
		&&\quad i=1,2,
		\quad\text{in $\QT$,}
		\\
		\label{Ki}
		\K_i &= \widehat\K_i - \sigma_i \uu,
		&&\quad i=1,2,
		\quad\text{on $\ST$,}
	\end{alignat}
	which represent the mass fluxes relative to the volume averaged velocity field. The mass balance equations \eqref{BL:MASS:B'} and \eqref{BL:MASS:S'} can thus be rewritten as
	\begin{alignat}{2}
		\label{BL:MASS:B*}
		\delt \rho_i + \Div(\rho_i\uu) + \Div \, \J_i &= 0,
		&&\quad i=1,2,
		\quad\text{in $\QT$,}
		\\
		\label{BL:MASS:S*}
		\deltn \sigma_i + \Divg(\sigma_i\uu) + \Divg \, \K_i  &= \widehat \J_{\Gamma,i} \cdot \n,
		&&\quad i=1,2,
		\quad\text{on $\ST$.}
	\end{alignat}
	In this context, $\J_i$ and $\K_i$, $i=1,2$, can be regarded as diffusive flow rates. We now define the total mass densities as
	\begin{alignat}{2}
		\label{DEF:RHO}
		\rho &= \rho_1 + \rho_2
		\quad\text{in $\QT$,}\\
		\label{DEF:SIGMA}
		\sigma &= \sigma_1 + \sigma_2
		\quad\text{in $\ST$.}
	\end{alignat}
	This leads to the relations
	\begin{alignat}{2}
		\label{BL:MASS:B**}
		\delt \rho + \Div(\rho\uu) + \Div(\J_1+\J_2)&= 0
		&&\quad\text{in $\QT$,}
		\\
		\label{BL:MASS:S**}
		\deltc \sigma + \sigma \Divg\uu  + \Divg(\K_1+\K_2) &= (\widehat \J_{\Gamma,1}  + \widehat \J_{\Gamma,2} ) \cdot \n
		&&\quad\text{on $\ST$.}
	\end{alignat}
	Moreover, in view of the balance of linear momentum with respect to the velocity field $\uu$, we have
	\begin{align}
		\label{BL:CLM}
		\delt (\rho\uu) + \Div(\rho\uu\otimes\uu) = \Div \, \tT 
		\quad\text{in $\QT$.}
	\end{align}
	Here, $\tT$ denotes the stress tensor, which needs to be specified by constitutive assumptions. In the following, we write 
	\begin{align}
		\label{DEF:JK}
		\J = \J_1 + \J_2 \quad\text{in $\QT$,}
		\qquad
		\K = \K_1 + \K_2 \quad\text{on $\ST$.}
	\end{align}
	Using \eqref{BL:MASS:B**} to reformulate \eqref{BL:CLM}, we obtain 
	\begin{align}
		\label{BL:CLM*}
		\begin{aligned}
			\rho(\delt\uu + \uu\cdot\Grad\uu) 
			&= \Div \, \tT + ( \Div \, \J )\uu
			= \Div(\tT + \uu\otimes \J) - \Grad\uu \; \J 
			\quad\text{in $\QT$.}
		\end{aligned}
	\end{align}
	Now, we define
	\begin{align}
		\label{DEF:TT}
		\T = \tT + \uu \otimes \J
		\quad\text{in $\QT$,}
	\end{align}
	Hence, \eqref{BL:CLM*} can be rewritten as
	\begin{align}
		\label{BL:CLM**}
		\rho \delt\uu + \Grad\uu \; (\rho\uu + \J) 
		= \Div\, \T 
		\quad\text{in $\QT$.}
	\end{align}
	We now multiply the equations \eqref{BL:MASS:B'} each by ${1}/{\trho_i}$, $i=1,2$, and we add the resulting equations. Recalling \eqref{DEF:PHIPSI*},\eqref{CONS:PHIPSI} and \eqref{DEF:VAVG:B}, we deduce 
	\begin{align}
		\label{PDE:DIVV}
		\begin{aligned}
			\Div \uu = \Div \left(\frac{\rho_1}{\trho_1}\uu_1 + \frac{\rho_2}{\trho_2}\uu_2 \right) = - \delt \left(\frac{\rho_1}{\trho_1} + \frac{\rho_2}{\trho_2}\right)
			= -\delt 1 = 0
			&&\quad\text{in $\QT$,}
		\end{aligned}
	\end{align}
	which means that the volume averaged velocity field is divergence-free. 
	Similarly, multiplying \eqref{BL:MASS:S'} by ${1}/{\tsigma_i}$, $i=1,2$, adding the resulting equations, and recalling \eqref{CONS:VOL}, \eqref{DEF:PHIPSI*}, \eqref{CONS:PHIPSI} and \eqref{DEF:VAVG:S}, we infer 
	\begin{align}
		\label{PDE:DIVV:S}
		\begin{aligned}
			& \Divg \, \uu 
			= \Divg \left(\frac{\sigma_1}{\tsigma_1}\uu_1 + \frac{\sigma_2}{\tsigma_2}\uu_2 \right)
			\\
			&\quad   
			= - \deltc \left(\frac{\sigma_1}{\tsigma_1} + \frac{\sigma_2}{\tsigma_2}\right) 
			+ \uu_\tau \cdot \Gradg \left(\frac{\sigma_1}{\tsigma_1} + \frac{\sigma_2}{\tsigma_2}\right)
			+ \left(\frac{\widehat\J_{\Gamma,1}}{\tsigma_1} + \frac{\widehat\J_{\Gamma,2}}{\tsigma_2}\right)\cdot\n
			\\
			&\quad  
			= \left(\frac{\widehat\J_{\Gamma,1}}{\tsigma_1} + \frac{\widehat\J_{\Gamma,2}}{\tsigma_2}\right)\cdot\n 
			= 0 
			&&\quad\text{on $\ST$.}
		\end{aligned}
	\end{align}
	We now define 
	\begin{alignat}{4}
		\label{DEF:TJ}
		\tJ_i &= \frac{\J_i}{\trho_i}, 
		&&\quad i=1,2,
		&\qquad \J_\phi &= \tJ_2-\tJ_1
		&&\quad\text{in $\QT$,}
		\\[1ex]
		\label{DEF:TJG}
		\tJ_{\Gamma,i} &= \frac{\J_{\Gamma,i}}{\tsigma_i}, 
		&&\quad i=1,2,
		&\qquad \J_\psi &= \tJ_{\Gamma,2}-\tJ_{\Gamma,1}
		&&\quad\text{on $\ST$,}
		\\[1ex]
		\label{DEF:TK}
		\tK_i &= \frac{\K_i}{\tilde{\sigma}_i},
		&&\quad i=1,2,
		&\qquad \K_\psi &= \tK_2-\tK_1
		&&\quad\text{on $\ST$.}
	\end{alignat}
	Multiplying the equations \eqref{BL:MASS:B*} each by ${1}/{\trho_i}$, $i=1,2$, and subtracting the resulting equations, and recalling \eqref{PDE:DIVV}, we deduce that
	\begin{align}
		\label{EQ:PHI}
        \deltb \phi + \Div \, \J_\phi
		= \delt \phi + \Div(\phi\uu) + \Div \, \J_\phi = 0
		\quad\text{in $\QT$.}
	\end{align}
	Moreover, multiplying the equations \eqref{BL:MASS:B*} each by ${1}/{\trho_i}$, adding the resulting equations, and recalling \eqref{CONS:PHIPSI} and \eqref{PDE:DIVV}, we conclude that
	\begin{align}
		\label{ID:DIVJ}
		\Div(\tJ_1 + \tJ_2) 
		= - \delt(\phi_1+\phi_2) - \Div\big((\phi_1+\phi_2)\uu \big) 
		= 0
		\quad\text{in $\QT$.}
	\end{align}
	Proceeding similarly with \eqref{BL:MASS:S*}, and recalling \eqref{CONS:PHIPSI} and \eqref{PDE:DIVV:S}, we obtain
	\begin{align}
		\label{EQ:PSI}
		\deltc \psi + \Divg \, \K_\psi
		= \deltn \psi  + \Divg(\psi\uu) + \Divg \, \K_\psi 
		= \J_\psi \cdot\n
		\quad\text{on $\ST$,}
	\end{align}
	and
	\begin{align}
		\label{ID:DIVK}
		\Divg(\tK_1 + \tK_2) 
		&= - \deltn (\psi_1+\psi_2) - \Divg\big((\psi_1+\psi_2)\uu \big) 
		+ \big(\tJ_{\Gamma,1} + \tJ_{\Gamma,2} \big)\cdot\n 
		= 0
	\end{align}
	on $\ST$.
	Using \eqref{DEF:TJ} and \eqref{ID:DIVJ}, which entails $\Div(\tJ_1)=-\Div(\tJ_2)$, we infer that
	\begin{align}
		\label{ID:DIVJ**}
		\begin{aligned}
			\Div \, \J 
			&= \frac 12 \bigg[ \trho_2 \Div \, \tJ_2 
			+ \trho_1 \Div \, \tJ_1 \bigg]
			+ \frac 12\bigg[ \trho_2 \Div \, \tJ_2 
			+ \trho_1 \Div \, \tJ_1 \bigg]
			\\[1ex]
			&= \frac{\trho_2-\trho_1}{2} \Div \, \tJ_2
			- \frac{\trho_2-\trho_1}{2} \Div \, \tJ_1
			= \frac{\trho_2-\trho_1}{2} \Div \, \J_\phi
			&&\quad\text{in $\QT$.}
		\end{aligned}
	\end{align}
	This means that there exists a divergence-free function $\J_0$ such that
    \begin{align}
        \label{ID:J:PRE}
		\J = \frac{\trho_2-\trho_1}{2} \J_\phi + \J_0
		\quad\text{in $\QT$.}
	\end{align}
    It is a common assumption in mixture theory that the relation between $\J$ and $\J_\phi$ is linear, that is, $\J_0 = \mathbf{0}$ in $\QT$. This means that 
	\begin{align}
		\label{ID:J}
		\J = \frac{\trho_2-\trho_1}{2} \J_\phi 
		\quad\text{in $\QT$.}
	\end{align}
    Alternatively, this linear relation can also be shown by assuming the final model to satisfy the global energy dissipation law \eqref{diss}. Using the relation \eqref{ID:J:PRE} would produce an additional summand $\Div(\uu \otimes \J_0)$ on the left-hand side of the Navier--Stokes equation \eqref{eqs:NSCH:1}. However, the energy dissipation law \eqref{diss} can only be satisfied if $\J_0$ is zero. Hence, in the following, we readily assume that $\J_0 = \mathbf{0}$ in $\QT$.
    
	Proceeding similarly to \eqref{ID:DIVJ**}, we further derive the identity
	\begin{align}
		\Divg \, \K = \frac{\tilde{\sigma}_2-\tilde{\sigma}_1}{2} \Divg \, \K_\psi
		\quad\text{on $\ST$.}
	\end{align}
	Assuming that $\phi$ and $\psi$ only attain values in the interval $[-1,1]$, we use \eqref{CONS:PHIPSI}, \eqref{DEF:PHIPSI} and \eqref{DEF:RHO} to conclude that the density $\rho = \rho(\phi)$ is given by
	\begin{align}
		\label{DEF:RHO*}
		\rho(\phi) = \frac{\trho_2-\trho_1}{2} \phi + \frac{\trho_2+\trho_1}{2}
		\quad\text{in $\QT$.}
	\end{align}

    \subsection{Option A: Model derivation via the Lagrange Multiplier Approach}
    \label{SUBSECT:LAG}

    One possibility to complete the model derivation is the \textit{Lagrange Multiplier Approach} (see~\cite{Liu2002}). It relies on local energy dissipation laws, the introduction of Lagrange multipliers, and constitutive assumptions to identify the unkown flux terms.
	
	\subsubsection{Local energy dissipation laws}
	
	As in \cite{Abels2012} and \cite{Giorgini2023}, we consider the following energy density in the bulk: 
	\begin{align}
		\label{DEF:END:B}
		e_{\Omega}(\uu,\phi,\Grad\phi) 
		= \frac{\rho(\phi)}{2}\abs{\uu}^2
		+ f(\phi,\Grad\phi)
		\quad\text{in $\QT$.}
	\end{align}
	Here, the first summand on the right-hand side represents the kinetic energy density. The second summand denotes the free energy density in the bulk. It is assumed to be of Ginzburg--Landau type, that is
	\begin{align}
		\label{DEF:GL:B}
		f(\phi,\Grad\phi)
		= \frac{\eps}{2}\abs{\Grad\phi}^2 + \frac 1\eps F(\phi)
		\quad\text{in $\QT$.}
	\end{align}
	As in \cite{Giorgini2023}, we additionally introduce the surface energy density
	\begin{align}
		\label{DEF:END:S}
		e_{\Gamma}(\phi,\psi,\Gradg\psi) 
		= g(\phi,\psi,\Gradg\psi)
		\quad\text{on $\ST$.}
	\end{align}
	Here, $g$ stands for the free energy density on the surface. It is also assumed to be of Ginzburg--Landau type with an additional term accounting for the relation between $\phi$ and $\psi$, namely
	\begin{align}
		\label{DEF:GL:S}
		g(\phi,\psi,\Gradg \psi) = \frac{\delta}{2}\abs{\Gradg\psi}^2
		+ \frac 1\delta G(\psi)
		+ \frac{1}{2} \h(K) (\alpha\psi-\phi)^2,
		\quad\text{$K\in[0,\infty]$}
	\end{align}
	on $\ST$, where $\h$ is the function introduced in \eqref{DEF:H}.
	Here, $\delta>0$ is related to the thickness of the diffuse interface on the boundary, 
	and $G$ is a potential that usually exhibits a double-well structure.
	
	To derive local energy dissipation laws, we now consider an arbitrary (sufficiently smooth) test volume $V(t)\subset \Om$,
	$t\in[0,T]$, that is transported by the velocity field $\uu$. 
	In an isothermal situation, the second law of thermodynamics leads to the dissipation inequality
	\begin{align}
		\label{IEQ:DISS}
		\begin{aligned}
			0 &\ge \ddt \left[ \int_{V(t)} e_{\Omega}(\uu,\phi,\Grad\phi) \dx
			+ \int_{\del V(t)\cap \Ga} e_{\Gamma}(\phi,\psi,\Gradg\psi) \dHH{d-1} \right]
			\\
			&\qquad + \int_{\del V(t)\cap \Om} \J_e \cdot \n_{\del V(t)} \dHH{d-1}
			+ \int_{\del (\del V(t)\cap \Ga)} \K_e \cdot \n_{\del (\del V(t)\cap \Ga)} \dHH{d-2}.
		\end{aligned}
	\end{align}
	In this inequality, the vector fields $\J_e$ and $\K_e$ are energy fluxes that will be specified later, $\n_{\del V(t)}$ is the outer unit normal vector field on $\del V(t)$, and $\n_{\del (\del V(t)\cap \Ga)}$ is the unit conormal vector field on the relative boundary $\del \big(\del V(t)\cap \Ga\big)$ within the submanifold $\Ga$. As we consider a thermodynamically closed system, there is no exchange of energy over the boundary $\Ga$ and thus, the domain of the first integral in the second line is just $\del V(t)\cap \Om$ instead of $\del V(t)$. Applying Gau\ss's divergence theorem for not necessarily tangential vector fields (see, e.g., \cite[Proposition~2.5.1]{BDGP-book}) on both integrals in the second line, we reformulate \eqref{IEQ:DISS} as
	\begin{align}
		\label{IEQ:DISS*}
		\begin{aligned}
			0 &\ge \ddt \left[ \int_{V(t)} e_{\Omega}(\uu,\phi,\Grad\phi) \dx
			+ \int_{\del V(t)\cap \Ga} e_{\Gamma}(\phi,\psi,\Gradg\psi) \dHH{d-1} \right]
			\\
			&\qquad + \int_{V(t)} \Div \, \J_e \dx
			+ \int_{\del V(t)\cap \Ga} \Divg \, \K_e 
			+ H \K_e \cdot \n
			- \J_e \cdot \n 
			\dHH{d-1}\,.
		\end{aligned}
	\end{align}
	Here, $H = -\frac{1}{d}\Divg\n$ is the mean curvature of $\Gamma(t)$, and we point out that $\n_{\del V(t)} = \n(t)$ on $\Gamma(t)$.
	Applying the Reynolds transport theorem for evolving domains (see, e.g., \cite[Theorem~2.11.1]{BDGP-book}), we deduce the identity
	\begin{align}
		\label{TRANSP:B}
		\ddt \int_{V(t)} e_{\Omega}(\uu,\phi,\Grad\phi) \dx
		= \int_{V(t)} \delt e_{\Omega}(\uu,\phi,\Grad\phi) 
		+ \Div\big(e_{\Omega}(\uu,\phi,\Grad\phi) \uu \big) \dx.
	\end{align}
	Similarly, recalling that $\Divg\uu=0$ (see~\eqref{PDE:DIVV:S}), we use the transport theorem for orientable evolving hypersurfaces (see, e.g., \cite[Theorem~2.10.1]{BDGP-book}) to obtain
	\begin{align}
		\label{TRANSP:S}
		\begin{split}
			&\ddt \int_{\del V(t)\cap \Ga} e_{\Gamma}(\phi,\psi,\Gradg\psi) \dHH{d-1} 
			= \int_{\del V(t)\cap \Ga} \deltc e_{\Gamma}(\phi,\psi,\Gradg\psi) 
			\dHH{d-1}.
		\end{split}
	\end{align}
	Now, by means of \eqref{TRANSP:B} and \eqref{TRANSP:S}, \eqref{IEQ:DISS*} can be reformulated as
	\begin{align}
		\label{IEQ:DISS**}
		\begin{aligned}
			0 &\ge \int_{V(t)} 
			\delt e_{\Omega}(\uu,\phi,\Grad\phi) 
			+ \Div\big(e_{\Omega}(\uu,\phi,\Grad\phi) \uu \big)
			+ \Div \, \J_e \dx
			\\
			&\quad 
			+ \int_{\del V(t)\cap \Ga} 
			\deltc e_{\Gamma}(\phi,\psi,\Gradg\psi) 
			\dHH{d-1}
			\\
			&\quad
			+ \int_{\del V(t)\cap \Ga} 
			\Divg \, \K_e 
			+ H \K_e \cdot \n 
			- \J_e\cdot\n 
			\dHH{d-1}.
		\end{aligned}
	\end{align}
	As \eqref{IEQ:DISS**} holds in particular for all test volumes $V(t)\subset \Om$ with $\del V(t) \cap \Ga = \emptyset$, we infer the following local energy dissipation law in the bulk:
	\begin{align}
		\label{IEQ:DISS:LOC:B}
		0 \ge -\mathcal D_{\Omega} 
		\coloneqq \delt e_{\Omega}(\uu,\phi,\Grad\phi) 
		+ \Div\big(e_{\Omega}(\uu,\phi,\Grad\phi) \uu \big)
		+ \Div \, \J_e 
        \quad\text{in $\QT$.}
	\end{align}
	Let now $c>0$ be arbitrary and let us consider a generic test volume $V(t)\subset \Om$ with $\abs{V(t)}$ being sufficiently small such that 
	\begin{align*}
		\int_{V(t)} \mathcal D_{\Omega} \dx < c.
	\end{align*}
	Invoking \eqref{IEQ:DISS:LOC:B}, we infer from \eqref{IEQ:DISS**} that
	\begin{align*}
		\int_{\del V(t)\cap \Ga} 
		\deltc e_{\Gamma}(\phi,\psi,\Gradg\psi) 
		+ \Divg \, \K_e 
		+ H \K_e \cdot \n 
		- \J_e\cdot\n \dHH{d-1}
		<c.
	\end{align*}
	As $c>0$ and the test volume $V(t)$ were arbitrary (except for the above smallness assumption), we conclude the following local energy dissipation law on the boundary:
	\begin{align}
		\label{IEQ:DISS:LOC:S}
		0 \ge -\mathcal D_{\Gamma}
		\coloneqq \deltc e_{\Gamma}(\phi,\psi,\Gradg\psi) 
		+ \Divg \, \K_e 
		+ H \K_e \cdot \n 
		- \J_e\cdot\n
        \quad\text{on $\ST$.}
	\end{align}
	
	\subsubsection{The Lagrange Multiplier Approach}\label{SUBSEC:LAG}
	
	Now, we introduce the Lagrange multipliers  $\mu$, $\theta$ and $q$  that need to be appropriately adjusted in the following approach. 
	In the final model, $\mu$ will represent the \textit{bulk chemical potential}, $\theta$ will represent the \textit{surface chemical potential} and $q$ will represent the \textit{surface pressure}.
	Invoking \eqref{PDE:DIVV}, \eqref{PDE:DIVV:S}, \eqref{EQ:PHI}, and \eqref{EQ:PSI}, the local energy dissipation laws \eqref{IEQ:DISS:LOC:B} and \eqref{IEQ:DISS:LOC:S} can be expressed as
	\begin{alignat}{2}
		\label{IEQ:DISS:B:0}
		0 \ge -\mathcal D_{\Omega} 
		&= \delt e_{\Omega}(\uu,\phi,\Grad\phi) 
		+ \Div\big(e_{\Omega}(\uu,\phi,\Grad\phi) \uu \big)
		+ \Div \, \J_e 
		\notag\\
		&\quad - \mu \underbrace{\big( 
		\delt \phi 
		+ \Div(\phi\uu) 
		+ \Div \, \J_\phi \big)}_{=0}
        - \, p \, \underbrace{\Div\uu \phantom{\big|\!} }_{=0}
		&&\quad\text{in $\QT$,}
		\\[1ex]
		\label{IEQ:DISS:S:0}
		0 \ge - \mathcal D_{\Gamma}
		&= \deltc e_{\Gamma}(\phi,\psi,\Gradg\psi) 
		+ \Divg \, \K_e 
		+ H \K_e \cdot \n 
		- \J_e\cdot\n
		\notag\\
		&\quad
		- \theta \underbrace{\big( 
		\deltc \psi 
		+ \Divg \, \K_\psi 
		- \J_\psi\cdot\n 
		\big)}_{=0}
		- \,q \, \underbrace{\Divg\uu \phantom{\big|\!} }_{=0}
		&&\quad\text{on $\ST$.}
	\end{alignat}%
	In the following, we will simply write $\rho$, $f$ and $g$ instead of $\rho(\phi)$, $f(\phi,\Grad\phi)$ and $g(\phi,\psi,\Gradg\psi)$ to provide a cleaner presentation.
	By the definition of the energy densities $e_{\Omega}$ and $e_{\Gamma}$ (see \eqref{DEF:END:B} and \eqref{DEF:END:S}), we reformulate \eqref{IEQ:DISS:B:0} and \eqref{IEQ:DISS:S:0} as
	\begin{alignat}{2}
		\label{IEQ:DISS:B:1}
		0 \ge -\mathcal D_{\Omega} 
		& = \delt \big(\tfrac12 \rho \abs{\uu}^2 \big)
		+ \Div\big(\tfrac12\rho \abs{\uu}^2 \uu \big)
		+ \delt f + \Div(f\uu)
		+ \Div \, \J_e
		\notag\\
		&\quad - \mu \big( 
		\delt \phi 
		+ \Div(\phi\uu) 
		+ \Div \, \J_\phi \big)
        + p \, \Div\uu
		&&\quad\text{in $\QT$,}
		\\[1ex]
		\label{EST:DISS:S:1}
		0 \ge - \mathcal D_{\Gamma}
		& = \deltc g 
		+ \Divg \, \K_e 
		+ H \K_e \cdot \n 
		- \J_e \cdot \n
		\notag\\
		&\quad- \theta \big( 
		\deltc \psi 
		+ \Divg \, \K_\psi 
		- \J_\psi\cdot\n 
		\big)
		+ q \, \Divg\uu
		&&\quad\text{on $\ST$.}
	\end{alignat}
	We now assume that $\psi$ is suitably extended onto a neighborhood of $\ST$.
	By means of the chain rule, the derivatives 
	$\delt f$ and $\deltc g$ 
	can be expressed as
	\begin{alignat}{2}
		\label{ID:DDTF}
		\delt f
		&= \delph f \, \delt \phi + \delgph f\, \delt \Grad\phi
		&&\quad\text{in $\QT$,}
		\\
		\label{ID:DDTGC}
		\deltc g
		&= \delph g \, \deltc \phi
		+ \delps g \, \deltc \psi 
		+ \delgps g\cdot \deltc \Gradg\psi
		&&\quad\text{on $\ST$.}
	\end{alignat}
	
	\paragraph{Computations in the bulk.}
	
	In the bulk, we proceed exactly as in \cite{Abels2012} and \cite{Giorgini2023} to reformulate inequality \eqref{IEQ:DISS:B:1} as
	\begin{align} 
		\label{EST:DISS:B:2}
		\begin{aligned}
			0 & \ge \Div\big[\J_e 
			- \tfrac 12 \abs{\uu}^2 \J 
			+ \T^T\uu 
			- \mu \J_\phi
			+ \delgph f \deltb \phi
			\big]
			\\
			& \quad + \big[ \delph f 
			- \Div(\delgph f) 
			- \mu 
			\big]
			\deltb \phi  
            + \Grad\mu \cdot \J_\phi
			\\
			& \quad - \big[ \T + p\I
			+ \Grad\phi \otimes \delgph f 
			\big] : \Grad\uu
			&&\quad\text{in $\QT$.}
		\end{aligned}
	\end{align}
	To ensure \eqref{EST:DISS:B:2}, we now choose the chemical potential $\mu$ and the energy flux $\J_e$ as
	\begin{alignat}{2}
		\label{DEF:MU*}
		\mu &= \delph f 
		- \Div\big( \delgph f \big)
		&&\quad\text{in $\QT$,}
		\\
		\label{DEF:JE}
		\J_e &= \frac{1}{2} \abs{\uu}^2 \J 
		- \T^T\uu 
		+ \mu \J_\phi
		- \delgph f \deltb \phi
		&&\quad\text{in $\QT$.}
	\end{alignat}
    Note that in \eqref{DEF:MU*}, we could potentially add a term $- \eta(\phi)\,\deltb\phi$ (with $\eta(\phi)>0$) on the right-hand side of \eqref{DEF:MU*} to account for dissipative friction. This would cause an additional dissipative term $-\eta(\phi)\abs{\deltb\phi} \le 0$ on the right-hand side of \eqref{EST:DISS:B:2}, which is consistent with the local energy dissipation law. In the final system of equations, this term would then lead to a viscous relaxation of the Cahn--Hilliard equation. However, in our model, we assume dissipative friction to be negligible and therefore, we do not include such a term.
	Due to the choices \eqref{DEF:MU*} and \eqref{DEF:JE}, the first two lines of the right-hand side of \eqref{EST:DISS:B:2} vanish.
	We further assume the mass flux $\J_\phi$ to follow Fick's law. This means that
	\begin{align}
		\label{DEF:JPHI}
		\J_\phi = -\mom(\phi) \Grad\mu \quad\text{in $\QT$,}
	\end{align}
	where $\mom = \mom(\phi)$ is a nonnegative function representing the mobility in the bulk.
	Hence, \eqref{EST:DISS:B:2} reduces to
	\begin{align}
		\label{EST:DISS:B:3}
		0  \ge  - \big[ \T + p\I
		+ \Grad\phi \otimes \delgph f 
		\big] : \Grad\uu
		- \mom(\phi) \abs{\Grad\mu}^2
		\quad\text{in $\QT$.}
	\end{align}
	We next assume that the stress tensor $\T$ can be expressed as
	\begin{align*}
		\T
        = \bS - p \I - \Grad\phi \otimes \delgph f 
		\quad\text{in $\QT$,}
	\end{align*}
	where the scalar variable $p$ denotes the pressure, $\I$ stands for the identity matrix, and $\mathbf S$ is the \textit{viscous stress tensor} that corresponds to irreversible changes of the energy due to friction. The term $\Grad\phi \otimes \delgph f$ can be interpreted as a Korteweg force. For Newtonian fluids, $\mathbf S$ is usually assumed to be given by
	\begin{align}
        \label{DEF:S}
		\mathbf S = 2 \nu(\phi) \D\uu,
	\end{align}
	where the nonnegative function $\nu=\nu(\phi)$ represents the viscosity of fluids, and $\D\uu$ denotes the symmetric gradient of $\uu$. This choice along with \eqref{DEF:MU*}, \eqref{DEF:JE} and the identity $p\I:\Grad \uu = p\, \Div \, \uu  = 0$ in $\QT$ ensures that \eqref{EST:DISS:B:2} is fulfilled as it reduces to
	\begin{align}
		\label{EST:DISS:B:4}
		0  \ge  - \nu(\phi) \abs{\D\uu}^2
		- \mom(\phi) \abs{\Grad\mu}^2
		\quad\text{in $\QT$.}
	\end{align}
	In particular, recalling \eqref{DEF:GL:B}, the chemical potential $\mu$ given by \eqref{DEF:MU*} can be expressed as
	\begin{align}
		\label{PDE:MU}
		\mu =  - \eps \Lap\phi + \frac 1\eps F'(\phi)
		\quad\text{in $\QT$.}
	\end{align}
	Moreover, the total stress tensor $\T$ is given by
	\begin{align}
		\label{DEF:T}
		\begin{aligned}
			\T 
			&= \bS 
			- p \I 
			- \eps \Grad\phi \otimes \Grad\phi
            \\
            &= 2 \nu(\phi) \D\uu 
			- p \I 
			- \eps \Grad\phi \otimes \Grad\phi
			\quad\text{in $\QT$.}
		\end{aligned}
	\end{align}
	Plugging \eqref{DEF:JPHI} into \eqref{EQ:PHI}, we further obtain the equation
	\begin{align}
		\label{PDE:PHI}
		\delt\phi + \Div(\phi\uu) = \Div\big( \mom(\phi) \Grad\mu \big)
		\quad\text{in $\QT$.}
	\end{align}
	Eventually, recalling \eqref{DEF:TT} and \eqref{DEF:T}, we use \eqref{BL:CLM} to derive the equation
	\begin{align}
		\label{PDE:NS}
		\delt (\rho\uu) 
		+ \Div\big(\uu\otimes(\rho(\phi)\uu + \J )\big) 
		- \Div\big(2\nu(\phi)\D\uu\big)
		+ \Grad p
		= - \eps\Div(\Grad\phi\otimes\Grad\phi)
		\quad\text{in $\QT$,}
	\end{align}
	where, according to \eqref{ID:J} and \eqref{DEF:JPHI}, we have
	\begin{align}
		\label{PDE:J}
		\J = - \frac{\trho_2-\trho_1}{2} \mom(\phi) \Grad\mu 
		\quad\text{in $\QT$.}
	\end{align}
	
	\paragraph{Computations on the boundary.}
	
	We now consider the local energy dissipation law \eqref{EST:DISS:S:1} on the surface.
	Recalling the formulae for $f$ (see~\eqref{DEF:GL:B}) and $\J_e$ (see~\eqref{DEF:JE}), we deduce from \eqref{EST:DISS:S:1} that
	\begin{align} \label{EST:DISS:S:2*}
		\begin{aligned}
			0 & \ge \delph g \, \deltc \phi
			+ \delps g \, \deltc \psi 
			+ \delgps g\cdot \deltc \Gradg\psi  
			+ \Divg \, \K_e
			+ H \K_e \cdot \n
			\\
			&\quad    
			- \tfrac 12 (\J\cdot\n)\uu\cdot\uu 
			- \T\n \cdot \uu
			- \mu \J_\phi\cdot \n
			+ \eps\, \deln\phi \, \deltc \phi    
			\\
			&\quad 
			- \theta \deltc \psi 
			- \Divg(\K_\psi \theta)
			+ \Gradg \theta\, \K_\psi 
			+ \theta \J_\psi\cdot \n
			\\
			&\quad 
			- \Divg(q \uu) + \Gradg q \cdot \uu 
		\end{aligned}
	\end{align}
	on $\ST$. 
    Now, the most technical step is to reformulate the term $\delgps g\cdot \deltc \Gradg\psi$. 
    According to \cite[Lemma~38]{BGN-book}, we have
    \begin{align} 
        \begin{split}
        \deltc \Gradg \psi - \Gradg \deltc \psi   
        &= \Gradg \uu \Gradg \psi - \mathbf{P}_\Gamma \big[ \Gradg \uu + (\Gradg \uu)^T \big] \mathbf{P}_\Gamma  \Gradg \psi
        \\
        &= \Gradg \uu  \Gradg \psi - \mathbf{P}_\Gamma \big[ \Gradg \uu + (\Gradg \uu)^T \big] \Gradg \psi
        \end{split}
    \end{align} 
    on $\ST$. Here, the symmetric matrix $\mathbf{P}_\Gamma$ represents the orthogonal projection onto the tangent space of $\Gamma=\Gamma(t)$.
    Invoking the definition of $g$ (see~\eqref{DEF:GL:S}), we further deduce
	\begin{align}
        \begin{aligned}
			&\delgps g \cdot \Big\{ \Gradg \uu  \Gradg \psi - \mathbf{P}_\Gamma \big[ \Gradg \uu + (\Gradg \uu)^T \big] \Gradg \psi \Big\}
            \\
            &\quad = -\delta \Gradg\psi \cdot (\Gradg \uu)^T \Gradg \psi 
            = - \delta \big( \Gradg\psi \otimes \Gradg\psi \big) : \Gradg \uu
		\end{aligned}
	\end{align}
    on $\ST$. From these two identities, we thus infer
    \begin{align}
        \begin{aligned}
			&\delgps g\cdot \deltc \Gradg\psi
            \\
            &\quad = \delgps g\, \Gradg \deltc \psi
                + \delgps g \cdot \Big\{ \Gradg \uu  \Gradg \psi - \mathbf{P}_\Gamma \big[ \Gradg \uu + (\Gradg \uu)^T \big] \Gradg \psi \Big\}
			\\
			&\quad = \delgps g\, \Gradg \deltc \psi
                - \delta \big( \Gradg\psi \otimes \Gradg\psi \big) : \Gradg \uu 
            \\[1ex]
            &\quad = \Divg\big[ \delgps g \deltc \psi - \delta (\Gradg\psi \otimes \Gradg\psi) \uu \big] 
                - \Divg\big( \delgps g \big) \deltc \psi 
                + \delta\, \Divg(\Gradg\psi \otimes \Gradg\psi) \uu
		\end{aligned}
	\end{align}
	on $\ST$. We further notice that
    \begin{align}
        \delta (\Gradg\psi \otimes \Gradg\psi) \uu \cdot \n
        = \delta (\Gradg\psi\cdot\uu) (\Gradg\psi\cdot\n)
        = 0
    \end{align}
    on $\ST$.
	Using these equations, inequality \eqref{EST:DISS:S:2*} can alternatively be expressed as
	\begin{align} \label{EST:DISS:S:2}
		\begin{aligned}
			0 & \ge \Divg\big[ \K_e
			- \theta \K_\psi 
            + \delgps g\, \deltc \psi 
            - \delta(\Gradg\psi \otimes \Gradg\psi) \uu
			- q \uu 
			\big]
			\\
			&\quad 
			+ \big[
			\delph g + \eps \Grad\phi\cdot \n
			\big] \deltc \phi
			+ \big[
			\delps g 
			- \Divg\big(\delgps g\big) 
			- \theta
			\big] \deltc \psi
            \\
			&\quad + \big[
			\T\n
			+ \Gradg q 
			- \tfrac 12 (\J\cdot\n)\uu 
            + \delta\, \Divg(\Gradg \psi \otimes \Gradg\psi) 
			\big]_\tau \cdot \uu_\tau
			\\
			&\quad 
			+ \big[
			\T\n\cdot\n 
			- \tfrac 12 (\J\cdot\n) (\uu\cdot\n) 
			\big] (\uu\cdot\n)
			\\
			&\quad
			+ H \K_e \cdot \n
			+ \Gradg \theta\, \K_\psi 
			+ \theta  \J_\psi \cdot \n - \mu \J_\phi\cdot \n
		\end{aligned}
	\end{align}
    on $\ST$.
	In order to ensure that the inequality \eqref{EST:DISS:S:2} is fulfilled, we choose the  mass flux $\K_\psi$, and the energy flux $\K_e$ as follows:
	\begin{alignat}{2}
		\label{DEF:KPSI}
		\K_\psi &= - \mga(\psi) \Gradg \theta
		&&\quad\text{on $\ST$,}
		\\
		\label{DEF:KE}
		\K_e &= \theta \K_\psi
        - \delgps g \, \deltc \psi
        + \delta(\Gradg\psi \otimes \Gradg\psi) \uu 
		+ q \uu
		&&\quad\text{on $\ST$.}
	\end{alignat}
	In \eqref{DEF:KPSI}, $\mga=\mga(\psi)$ is a nonnegative function representing the mobility of the materials on the surface.
	The choice of $\K_e$ entails that the first line of the right-hand side in \eqref{EST:DISS:S:2} vanishes. Except for $q\uu$, all terms on the right-hand side of \eqref{DEF:KE} are tangential vector fields. We thus have
	\begin{align}
		\label{EQ:KE}
        H\K_e \cdot \n 
        = Hq (\uu\cdot\n)
         \quad\text{on $\ST$.}   
	\end{align}
	Due to \eqref{DEF:KPSI}--\eqref{EQ:KE}, inequality \eqref{EST:DISS:S:2} now reduces to
	\begin{align} \label{EST:DISS:S:3*}
        \begin{aligned}
			0 & \ge  \big[
			    \delph g + \eps \Grad\phi\cdot \n
			\big] \deltc \phi
			+ \big[
    			\delps g 
    			- \Divg\big(\delgps g\big) 
    			- \theta
			\big] \deltc \psi
            \\
			&\quad + \big[
    			\T\n
    			+ \Gradg q 
    			- \tfrac 12 (\J\cdot\n)\uu 
                + \delta\, \Divg(\Gradg \psi \otimes \Gradg\psi) 
			\big]_\tau \cdot \uu_\tau
			\\
			&\quad 
			+ \big[
    			\T\n\cdot\n 
    			- \tfrac 12 (\J\cdot\n) (\uu\cdot\n) 
                + Hq
			\big] (\uu\cdot\n)
			\\
			&\quad
			- \mga(\psi)\abs{\Gradg\theta}^2 
			+ \theta  \J_\psi \cdot \n - \mu \J_\phi\cdot \n
		\end{aligned}
	\end{align}
	on $\ST$. We now assume that the flux terms $\J_\phi$ and $\J_\psi$ are directly proportional. This means that there exists a coefficient $\beta$ such that
	\begin{equation}
        \label{ASS:FLUX}
		\J_\psi \cdot \n = \beta \J_\phi \cdot \n = - \beta \, m_\Omega(\phi) \Grad\mu \cdot \n \quad \ton \ST.
	\end{equation}
	Here, the second equality follows from \eqref{DEF:JPHI}.
    The idea behind this assumption is that $\beta$ 
    represents the fraction of bulk material that is effectively transformed into surface material when reaching the boundary. In this paper, we assume for simplicity that $\beta$ is a constant. However, in general, the coefficient could also be a sufficiently regular function $\beta:Q_T\to\R$ (see, e.g., \cite{Goldstein2011}).
    
	We further make the constitutive assumptions
	\begin{alignat}{2}
		\label{BC:NAV}
		\big[
		      \T\n
			+ \Gradg q 
			- \tfrac 12 (\J\cdot\n)\uu 
            + \delta\, \Divg(\Gradg \psi \otimes \Gradg\psi) 
		\big]_\tau 
		&=
		-\gamma_\tau\uu_\tau
		&&\ton \ST,
		\\
		\label{BC:NOR}
		      \T\n\cdot\n 
            - \tfrac 12 (\J\cdot\n) (\uu\cdot\n) 
            + Hq
		&=
		-\gamma_\n (\uu\cdot\n)
		&&\ton \ST,
	\end{alignat}
	where $\gamma_\tau$ and $\gamma_\n$ are nonnegative functions, which may depend on the phase fields. Equation \eqref{BC:NAV} can be interpreted as a \textit{generalized Navier slip boundary condition} with $\gamma_\tau$ being the slip parameter. We point out that an analogous boundary condition appeared in \cite{Giorgini2023}.
    The boundary condition \eqref{BC:NOR} relates the normal component of $\uu$ to the normal traction $\T\n$ and the geometry of the evolving boundary through its mean curvature $H$. The coefficient $\gamma_\n>0$ acts as a weighting parameter for the energy dissipation due to boundary deformation (cf.~\eqref{diss}).
	By means of \eqref{ASS:FLUX}--\eqref{BC:NOR}, inequality
	\eqref{EST:DISS:S:3*} further reduces to 
	\begin{align} \label{EST:DISS:S:3}
		\begin{aligned}
			0 &\ge 
			\big[
			\delph g + \eps \Grad\phi\cdot \n
			\big] \deltc \phi
			+ \big[
			\delps g 
			- \Divg\big(\delgps g\big) 
			- \theta
			\big] \deltc\psi
			\\
			&\quad
			- \gamma_\tau\abs{\uu_\tau}^2 
			- \gamma_\n (\uu\cdot\n)^2
			- \mga(\psi) \abs{\Gradg \theta}^2
			\\
			&\quad
			- (\beta\theta - \mu)\,  \mom(\phi) \Grad\mu \cdot \n
		\end{aligned}
	\end{align}
	on $\ST$.
	The terms in the second line are clearly nonpositive. In order to ensure that the third line is also nonpositive, we assume that one of the following boundary conditions holds:
	\begin{subequations} 
		\label{BC:MUTHETA}
		\begin{alignat}{2}
			\beta\theta - \mu &= 0
			&&\quad\text{on $\ST$,}
			\\
			\mom(\phi) \Grad\mu \cdot \n &= \tfrac 1L (\beta\theta - \mu)
			&&\quad\text{on $\ST$ for a constant $L\in(0,\infty)$,}
			\\
			\mom(\phi) \Grad\mu \cdot \n &= 0
			&&\quad\text{on $\ST$.}
		\end{alignat}
	\end{subequations}
	It thus remains to ensure that the first line in \eqref{BC:MUTHETA} is nonpositive. 
    The approach depends on the parameter $K$ that appears in the energy density $g$ (see \eqref{DEF:GL:S})
    and we need to handle the cases $K=0$ and $K\in(0,\infty]$ separately.
    \begin{subequations}
    \begin{itemize}[leftmargin=*]
    	\item\textit{The case $K=0$.}
    	In this case, we assume that $\phi$ and $\psi$ are directly proportional on $\ST$. More precisely, we assume that
        \label{BC:PHIPSI}
        \begin{align}
            \label{BC:PHIPSI:1}
            \alpha\psi-\phi = 0 \quad\ton\ST 
        \end{align}
        with a prescribed parameter $\alpha\in\R\setminus\{0\}$.
        In the definition of $g$, we have $\h(K) = \h(0) = 0$ and thus $\delph g = 0$ on $\ST$. Moreover, due to \eqref{BC:PHIPSI:1}, $\alpha^{-1}\phi$ can be interpreted as an extension of $\psi$ in some neighborhood of $\ST$. This entails $\alpha\Grad\psi = \Grad\phi$ on $\ST$ and we thus have
        \begin{equation*}
            \deltc\phi = \alpha\deltc\psi
            \quad\text{on $\ST$.}
        \end{equation*}
        Consequently, the first line of \eqref{EST:DISS:S:3} can be rewritten as
        \begin{align*}
            \big[
            \delps g 
            - \Divg\big(\delgps g\big) 
            - \theta
            + \alpha\eps \Grad\phi\cdot \n
            \big] \deltc\psi
            \quad\ton\ST.
        \end{align*}
        Choosing 
        \begin{align*}
            \theta 
            &= - \Divg\big(\delgps g\big)
            + \delps g 
            + \alpha\eps \Grad\phi\cdot \n,
        \end{align*}
        we thus ensure that the first line of \eqref{EST:DISS:S:3} vanishes, which shows that the inequality \eqref{EST:DISS:S:3} is fulfilled.
        
        \item \textit{The case $K\in(0,\infty]$.}
        In this case, we assume that $\phi$ and $\psi$ are coupled by the boundary condition
        \begin{alignat}{2}
            \label{BC:PHIPSI:2}
            K\eps \deln\phi &= \alpha\psi-\phi
            &&\quad\ton\ST \text{ if $K\in (0,\infty)$},
            \\
            \label{BC:PHIPSI:3}
            \deln\phi &= 0
            &&\quad\ton\ST \text{ if $K=\infty$},
        \end{alignat}
        where $\alpha\in\R$ is a prescribed parameter.
    	In view of the definition of $g$ (see \eqref{DEF:GL:S}), this ensures that 
    	\begin{alignat*}{2}
    		\delph g  
    		&= - \h(K) (\alpha\psi - \phi)  
    		= - \eps \Grad\phi\cdot\n 
    		&&\ton\ST,
    		\\
    		\delps g  
    		&= \delta^{-1} G'(\psi) + \alpha \h(K) (\alpha\psi - \phi)  
    		=  \delta^{-1} G'(\psi) + \alpha \eps \Grad\phi\cdot\n 
    		&&\ton\ST.
    	\end{alignat*}
    	Consequently, the first line of \eqref{EST:DISS:S:3} can be rewritten as
    	\begin{align*}
    		\big[
    		\delps g 
    		- \Divg\big(\delgps g\big) 
    		- \theta
    		\big] \deltc\psi
    		\quad\ton\ST.
    	\end{align*}
    	Choosing 
    	\begin{align*}
    		\theta 
    		&= - \Divg\big(\delgps g\big)
    		+ \delps g ,
    	\end{align*}
    	we thus ensure that the first line of \eqref{EST:DISS:S:3} vanishes, which shows that the inequality \eqref{EST:DISS:S:3} is fulfilled.
	\end{itemize}
    \end{subequations}
	By the definition of $g$ (see \eqref{DEF:GL:S}), we conclude in all cases $K\in[0,\infty]$ that $\theta$ is given by
	\begin{align}
		\label{PDE:THETA}
		\theta = - \delta \Lapg\psi + \frac 1\delta G'(\psi) + \alpha \eps \deln\phi 
		\ton\ST.
	\end{align}
	Eventually, by substituting \eqref{DEF:JPHI} and \eqref{DEF:KPSI} into \eqref{EQ:PSI}, we obtain
	\begin{align}
		\label{PDE:PSI}
		\deltc \psi + \Divg(\psi\uu_\tau) = \Divg(\mga(\psi) \Gradg \theta) - \beta \mom(\phi) \Grad\mu\cdot\n
		\quad\text{on $\ST$.}
	\end{align}

    In summary, we have thus derived system \eqref{eqs:NSCH} along with formulae \eqref{eq:rho} and \eqref{eq:T} with the following correspondences:
    \begin{center}
    \makeatletter \setlength{\tabcolsep}{2pt} \makeatother
    \begin{tabular}{rclcrclcrcl}
        \eqref{eqs:NSCH:1} 
        &$\widehat{=}$ 
        &\eqref{PDE:NS}, \eqref{PDE:DIVV}, 
        &$\phantom{xx}$ 
        &\eqref{eqs:NSCH:2} 
        &$\widehat{=}$ 
        &\eqref{PDE:PHI}, 
        &$\phantom{xx}$ 
        &\eqref{eqs:NSCH:3} 
        &$\widehat{=}$ 
        &\eqref{PDE:MU}, 
        \\
        \eqref{eqs:NSCH:4} 
        &$\widehat{=}$ 
        &\eqref{PDE:PSI},
        &$\phantom{xx}$ 
        &\eqref{eqs:NSCH:5} 
        &$\widehat{=}$ 
        &\eqref{PDE:THETA},
        &$\phantom{xx}$ 
        &\eqref{eqs:NSCH:6} 
        &$\widehat{=}$ 
        &\eqref{BC:PHIPSI}, \eqref{BC:MUTHETA}, 
        \\
        \eqref{eqs:NSCH:7} 
        &$\widehat{=}$ 
        &\eqref{BC:NAV}, 
        &$\phantom{xx}$ 
        &\eqref{eqs:NSCH:8} 
        &$\widehat{=}$ 
        &\eqref{BC:NOR},
        &$\phantom{xx}$ 
        &\eqref{eqs:NSCH:9} 
        &$\widehat{=}$ 
        &\eqref{PDE:DIVV:S}, \eqref{BC:NORVEL},
        \\
        \eqref{eq:rho} 
        &$\widehat{=}$ 
        &\eqref{DEF:RHO*},
        &$\phantom{xx}$ 
        &\eqref{eq:J} 
        &$\widehat{=}$ 
        &\eqref{PDE:J},
        &$\phantom{xx}$ 
        &\eqref{eq:T} 
        &$\widehat{=}$ 
        &\eqref{DEF:T}, \eqref{DEF:S}.
    \end{tabular}
    \end{center}
    Therefore, our model derivation is complete.

    
    \subsection{Option B: Model derivation via the Energetic Variational Approach}

    Another possibility to complete the model derivation is the \textit{Energetic Variational Approach (EnVarA)} (see, e.g., \cite{Hyon2010,Giga2018,Wang2022}). It relies on the balance of inertial, conservative and dissipative forces. The inertial and conservative forces are determined by means of the \textit{least action principle}. Assuming a global energy dissipation rate, the dissipative forces are derived via \textit{Onsager's principle of maximal energy dissipation}.
    The a priori unknown flux terms $\J_\phi$, $\K_\psi$ and $\J_\psi$, which appear in the balance laws \eqref{EQ:PHI} and \eqref{EQ:PSI}, are determined by means of the dissipation law for the total energy functional.

    The EnVarA strongly relies on considerations of the \textit{flow map} associated with the fluid mixture
    \begin{align*}
        \x:[0,T]\times \overline{\Omega(0)} \to \overline{\Omega(t)},
    \end{align*}
    which is defined as the unique solution to the initial value problem
    \begin{align*}
        \ddt \x(t,X) = \uu\big(t,\x(t,X)\big), \qquad
        \x(0,X) = X.
    \end{align*}
    It is a key assumption of our model that the bulk and the surface dynamics are governed by the same flow map that corresponds to the velocity field $\uu$.

    Compared to the Lagrange Multiplier Approach presented in Subsection~\ref{SUBSECT:LAG}, the EnVarA has the advantage that the inertial, conservative and dissipative forces can be identified clearly. 
    However, a drawback is that we are only able to rigorously perform the derivation by the EnVarA only under two additional assumptions:
    \begin{itemize}[leftmargin=*, topsep=0em, partopsep=0em, parsep=1ex, itemsep=0ex]
    
        \item \textbf{No mass flux between bulk and surface.}
        As in \cite{Liu2019}, we need to assume that no transfer of material between the bulk and the surface is allowed. This means we assume that
        \begin{align}
            \label{ASS:NOFLUX}
            \J_\phi \cdot \n = 0
            \quad\text{and}\quad
            \J_\psi \cdot \n = 0
            \ton\ST.
        \end{align}
        Consequently, the balance laws \eqref{EQ:PHI} and \eqref{EQ:PSI} can be restated as
        \begin{alignat}{2}
            \label{EQ:PHI:E*}
            \deltb \phi + \Div \, \J_\phi &= 0
            &&\tin\QT,
            \\
            \label{EQ:PSI:E*}
            \deltc \psi + \Divg \, \K_\psi &= 0
            &&\ton\ST,
        \end{alignat}
        which will lead to the parameter choice $L=\infty$ in the final model \eqref{eqs:NSCH}.
        We point out that the neglection of mass exchange between bulk and surface was a crucial assumption in \cite{Liu2019} to carry out the EnVarA. To the best of our knowledge, there is no contribution in the literature so far, where the EnVarA allows for the inclusion of a transfer of material between bulk and surface. 
        
        \item \textbf{Matched densities.} We further assume that both fluids have the same individual density, that is, $\trho_1 = \trho_2$. 
        In view of \eqref{ID:J} and \eqref{DEF:RHO*}, this assumption leads to 
        \begin{align}
            \label{SIMP:JRHO}
            \rho = \rho(\phi) \equiv \trho_1 = \trho_2
            \quad\text{and}\quad
            \J = \mathbf{0}
            \tin\QT.
        \end{align}
        The assumption of matched densities will allow us to compute the inertial forces by proceeding, for example, as in \cite{Jiang2017}. We point out that in \cite{Liu2015}, the authors derived a Navier--Stokes--Cahn--Hilliard model even in the case of unmatched densities. However, they were only able to recover the Navier--Stokes equation in its weak formulation, which is not enough to rigorously identify the inertial forces. It is also not completely clear whether their techniques are directly transferable to our situation.
        For these reasons, we restrict ourselves to only consider the situation of matched densities in this subsection.
	\end{itemize}

    \subsubsection{Energy dissipation}

    A key ingredient of the EnVarA is the assumption of a global energy dissipation law.
    To this end, we introduce the total energy functional
    \begin{align*}
        \cE(\uu,\phi,\psi)
        \coloneqq \cE^\mathrm{kin}(\uu,\phi) + \cE^\mathrm{free}(\phi,\psi),
    \end{align*}
    where the \textit{kinetic energy} $\cE^\mathrm{kin}$, the \textit{bulk free energy} $\cE^\mathrm{bulk}$, and the \textit{surface free energy}   $\cE^\mathrm{surf}$ are given by
    \begin{align*}
        \cE^\mathrm{kin}(\uu) 
        &\coloneqq \intO \frac{\rho}{2} \abs{\uu}^2 \dx\,,
        \\
        \cE^\mathrm{bulk}(\phi)
        &\coloneqq \intO f(\phi,\Grad\phi) \dx\,,
        \\
        \cE^\mathrm{surf}(\phi,\psi)
        &\coloneqq \intG g(\phi,\psi,\Gradg \psi) \dHH{d-1}\,.
    \end{align*}
    Here, the free energy densities $f$ and $g$ are defined as in \eqref{DEF:GL:B} and \eqref{DEF:GL:S}, respectively. 
    This definition is consistent with \eqref{DEF:E}. 
    As we consider the situation of matched densities, we simply write $\rho$ instead of $\rho(\phi)$, bearing in mind that $\rho$ is constant (cf.~\eqref{SIMP:JRHO}).
    We assume that the energy dissipation law
    \begin{align}
        \label{EN:DISS:GLOB}
        \ddt \cE = - \cD
    \end{align}
    holds, where the dissipation rate $\cD$ is given by
    \begin{align}
        \label{DEF:DISS:GLOB}
        \begin{split}
        \cD 
        &\coloneqq \intO 2 \nu(\phi) \abs{\D \uu}^2 \dx
			+ \intG \gamma_{\tau} \abs{\uu_{\tau}}^2 \dHH{d-1} 
			+ \intG \gamma_{\n} \abs{\uu \cdot \n}^2 \dHH{d-1}
        \\
        &\quad
            + \intO \frac{1}{m_{\Omega}(\phi)} \abs{\J_\phi}^2 \dx
			+ \intG \frac{1}{m_{\Gamma}(\psi)} \abs{\K_\psi}^2 \dHH{d-1}
        .
        \end{split}
    \end{align}
    Here, $\nu(\phi)$ denotes the \textit{viscosity} of the mixture, which may depend on the phase-field. The corresponding term in the dissipation rate accounts for internal friction within the fluids. The parameter $\gamma_\tau>0$ is related to diffusion via tangential friction of the materials at the boundary, and recalling $\uu\cdot\n = \cV$ on $S_T$ (see \eqref{BC:NORVEL}), $\gamma_\n>0$ acts as a weighting parameter for energy dissipation due to boundary deformation.
    
    The integrals involving $\J_\phi$ and $\K_\psi$ can be interpreted as a bulk-surface generalization of the usual quadratic bulk dissipation rate that is motivated by Onsager’s linear response theory \cite{Onsager1931-I,Onsager1931-II}. 
    Here, the positive  $m_{\Omega}=m_{\Omega}(\phi)$ and $m_{\Gamma}=m_{\Gamma}(\psi)$ represent the mobility in the bulk and on the surface, respectively. 

    \begin{remark}
        Alternatively, proceeding as in \cite{Liu2019} (see also, e.g., \cite{Xu2014,LiuPei2018,Metzger2020}), the flux terms $\J_\phi$ and $\K_\psi$ could be represented by means of microscopic effective velocities $\uu_\phi$ and $\uu_\psi$, respectively. 
        This means that
        \begin{alignat*}{2}
            \J_\phi &= \phi\, \uu_\phi
            \tin\QT,
            \\
            \K_\psi &= \psi\, \uu_\psi
            \ton\ST,
        \end{alignat*}
        In this way, \eqref{EQ:PHI:E*} and \eqref{EQ:PSI:E*} can be alternatively expressed as continuity equations, namely
        \begin{alignat}{2}
            \label{EQ:PHI:E:C}
            \delt \phi + \Div\big( \phi (\uu + \uu_\phi)\big) &= 0
            &&\tin\QT,
            \\
            \label{EQ:PSI:E:C}
            \delt \psi + \Divg\big( \psi (\uu + \uu_\psi) \big) &= 0
            &&\ton\ST.
        \end{alignat}
        In this case, assuming the dissipation rate to be given by
        \begin{align}
            \label{DEF:DISS:GLOB*}
            \begin{split}
            \cD 
            &= \intO 2 \nu(\phi) \abs{\D \uu}^2 \dx
    			+ \intG \gamma_{\tau} \abs{\uu_{\tau}}^2 \dHH{d-1} 
    			+ \intG \gamma_{\n} \abs{\uu \cdot \n}^2 \dHH{d-1}
            \\
            &\quad
                + \intO \phi^2 \frac{\uu_\phi}{m_{\Omega}(\phi)} \cdot \uu_\phi \dx
    			+ \intG \psi^2 \frac{\uu_\psi}{m_{\Gamma}(\psi)} \cdot \uu_\psi \dHH{d-1}
            \end{split}
        \end{align}
        (cf.~\cite{Liu2019}), the EnVarA leads to the same final model as in our flux-based formulation.
    \end{remark}

    \subsubsection{Derivation of the inertial and conservative forces via the least action principle}

    To derive the inertial and conservative forces, which correspond to reversible dynamics in our mechanical system, we use the \textit{least action principle}. First, we first introduce the \textit{Lagrangians} associated with the kinetic energy and the bulk and surface free energies. They read as follows:
    \begin{align*}
        L^\mathrm{kin}\big( \x(\cdot) \big)
        &\coloneqq \intO \frac{\rho}{2} \abs{\uu}^2 \dx \,,
        \\
        L^\mathrm{bulk}\big( \x(\cdot) \big)
        &\coloneqq - \intO f(\phi,\Grad\phi) \dx \,,
        \\
        L^\mathrm{surf}\big( \x(\cdot) \big)
        &\coloneqq - \intG g(\phi,\psi,\Gradg \psi) \dHH{d-1}\,.
    \end{align*}
    Here, the notation ``$\x(\cdot)$'' is used to express that the Lagrangians are to be understood as functions of the flow map $\x:[0,T]\times \overline{\Omega(0)} \to \overline{\Omega(t)}$.
    Next, we define the \textit{action functional}
    \begin{align*}
        \cA\big( \x(\cdot) \big)
        = \cA^\mathrm{kin}\big( \x(\cdot) \big) 
            + \cA^\mathrm{bulk}\big( \x(\cdot) \big) 
            + \cA^\mathrm{surf}\big( \x(\cdot) \big)\,,
    \end{align*}
    where
    \begin{alignat*}{2}
        \cA^\mathrm{kin}\big( \x(\cdot) \big)
        &\coloneqq \int_0^T L^\mathrm{kin}\big( \x(\cdot) \big) \dt
        &&= \int_0^T \intO \frac{\rho}{2} \abs{\uu}^2 \dx \dt\,,
        \\
        \cA^\mathrm{bulk}\big( \x(\cdot) \big)
        &\coloneqq \int_0^T L^\mathrm{bulk}\big( \x(\cdot) \big) \dt
        &&= - \int_0^T \intO f(\phi,\Grad\phi) \dx \dt\,,
        \\
        \cA^\mathrm{surf}\big( \x(\cdot) \big)
        &\coloneqq \int_0^T L^\mathrm{surf}\big( \x(\cdot) \big) \dt
        &&= - \int_0^T \intG g(\phi,\psi,\Gradg \psi) \dHH{d-1} \dt\,.
    \end{alignat*}
    
    The least action principle states that the first variation of the action functional is zero. To compute this first variation, let 
    \begin{align*}
        \y:[0,T]\times \overline{\Omega(0)} \to \R^d
    \end{align*}
    be an arbitrary function, which satisfies $\y(0,\cdot) = \y(T,\cdot) = \mathbf{0}$ in $\overline{\Omega(0)}$.
    Furthermore, we write $\tilde\y$ to denote its Eulerian counterpart, that is,
    \begin{alignat*}{2}
        \y(t,X) &= \tilde\y(t,\x(t,X)) 
        &&\quad\text{for all $(t,X)\in [0,T]\times \overline{\Omega(0)}$.}
    \end{alignat*}
    To account for the incompressibility constraints \eqref{PDE:DIVV} and \eqref{PDE:DIVV:S}, we further assume that
    \begin{alignat}{2}
        \label{COND:Y:B}
        \Div(\tilde \y) &= 0 
        &&\quad\ton\QT,
        \\
        \label{COND:Y:S}
        \Divg(\tilde \y) &= 0
        &&\quad\ton\ST.
    \end{alignat} 
    For any $\varsigma\in (-1,1)$, we set
    \begin{align*}
        \x^\varsigma(\cdot)\coloneqq \x(\cdot) + \varsigma\, \y(\cdot)\,.
    \end{align*}
    
    Proceeding as in \cite{Jiang2017}, the first variation of the kinetic part of the action functional with respect to the flow map can be computed as
    \begin{align*}
        &\delta_{\x(\cdot)} \cA^\mathrm{kin}\big( \x(\cdot) \big) \big[ \y(\cdot) \big]
        = \frac{\mathrm{d}}{\mathrm{d} \varsigma}\Big|_{\varsigma=0} \cA^\mathrm{kin}\big( \x^\varsigma(\cdot) \big)
        \\
        &\quad= \frac{\mathrm{d}}{\mathrm{d} \varsigma}\Big|_{\varsigma=0} \int_0^T \int_{\Omega(0)} 
            \frac 12 \rho\, \abs{\partial_t \x^\varsigma}^2 \dX\dt
        = \int_0^T \int_{\Omega(0)} 
            \rho\, \partial_t \x  \cdot \partial_t \y \dX\dt
        \\
        &\quad = - \int_0^T \int_{\Omega(0)} 
            \rho\, \partial_t \partial_t \x  \cdot \y \dX\dt
        = - \int_0^T \intO 
            \big[\rho \partial_t\uu + \rho (\uu \cdot \Grad)\uu \big]
            \cdot \tilde\y \dx\dt
        \\
        &\quad= - \int_0^T \intO 
            \big[ \rho \partial_t \uu + \Div( \uu \otimes \rho \uu )\big]
            \cdot \tilde\y \dx\dt\,.
    \end{align*}

    Next, 
    we can proceed analogously to \cite[Sect.~7, Step~1, Eq.~(7.3)]{Liu2019} to compute the first variation of $\cA^\mathrm{bulk}$. However, in our case, many terms actually vanish due to the incompressibility constraint \eqref{COND:Y:B}. We obtain
    \begin{align*}
        &\delta_{\x(\cdot)} \cA^\mathrm{bulk}\big( \x(\cdot) \big) \big[ \y(\cdot) \big]
        = \frac{\mathrm{d}}{\mathrm{d} \varsigma}\Big|_{\varsigma=0} \cA^\mathrm{bulk}\big( \x^\varsigma(\cdot)\big)
        \\
        &\quad= \int_0^T \intO \partial_{\Grad\phi} f \cdot \Big[ (\Grad \tilde\y)^T \, \Grad \phi \Big]\dx\dt
        = \int_0^T \intO \eps (\Grad\phi \otimes \Grad\phi) : \Grad \tilde\y \dx\dt 
        \\
        \begin{split}
        &\quad= - \int_0^T \intO \eps\, \Div(\Grad\phi \otimes \Grad\phi) \cdot \tilde\y \dx\dt
            + \int_0^T \intG 
            \Big[ \eps\, (\Grad\phi \otimes \Grad\phi)\n \Big]_\tau \cdot \tilde\y \dHH{d-1} \dt
        \\
        &\quad\qquad 
            + \int_0^T \intG 
            \Big[ \eps\, (\Grad\phi \otimes \Grad\phi)\n \cdot \n\Big] \n \cdot \tilde\y  \dHH{d-1} \dt\,.
        \end{split}
    \end{align*}

    Lastly, 
    we can proceed analogously to \cite[Sect.~7, Step~1, Eq.~(7.5)]{Liu2019} to compute the first variation of $\cA^\mathrm{surf}$. 
    As in the bulk, many terms actually drop out due to the incompressibility constraint \eqref{COND:Y:S}. We obtain
    \begin{align*}
        &\delta_{\x(\cdot)} \cA^\mathrm{surf}\big( \x(\cdot) \big) \big[ \y(\cdot) \big]
        = \frac{\mathrm{d}}{\mathrm{d} \varsigma}\Big|_{\varsigma=0} \cA^\mathrm{surf}\big( \x^\varsigma(\cdot) \big)
        \\
        &\quad= \int_0^T \intG \partial_{\Gradg\psi} g \cdot \Big[ (\Grad \tilde\y)^T \, \Gradg \psi \Big]\dHH{d-1}\dt
        \\
        &\quad= \int_0^T \intG \delta (\Gradg\psi \otimes \Gradg\psi) : \Gradg \tilde\y \dHH{d-1}\dt 
        \\
        &\quad= - \int_0^T \intG \Big[\delta\, \Divg(\Gradg\psi \otimes \Gradg\psi)\Big]_\tau \cdot \tilde\y \dHH{d-1}\dt\,.
    \end{align*}

    Hence, in summary, we have
    \begin{align*}
        &\delta_{\x(\cdot)} \cA\big( \x(\cdot) \big) \big[ \y(\cdot) \big]
        = \frac{\mathrm{d}}{\mathrm{d} \varsigma}\Big|_{\varsigma=0} \cA\big( \x^\varsigma(\cdot) \big)
        \\
        &\quad= 
            - \int_0^T \intO 
            \Big[ 
                \rho\, \big( \delt \uu + (\uu\cdot\Grad) \uu \big)
            \Big] \cdot \tilde\y \dx\dt
        \\
        &\quad\qquad 
            - \int_0^T \intO \eps\, \Div(\Grad\phi \otimes \Grad\phi) \cdot \tilde\y \dx\dt
        \\
        &\quad\qquad 
            + \int_0^T \intG 
            \Big[ 
                \eps\, (\Grad\phi \otimes \Grad\phi)\n 
                - \delta\, \Divg(\Gradg\psi \otimes \Gradg\psi)
            \Big]_\tau \cdot \tilde\y \dHH{d-1} \dt
        \\
        &\quad\qquad
            + \int_0^T \intG 
            \Big[ 
                \eps\, (\Grad\phi \otimes \Grad\phi)\n \cdot \n
            \Big] \n \cdot \tilde\y \dHH{d-1} \dt.
    \end{align*}
    Invoking the least action principle, which states that the first variation of the action functional is zero, we derive the inertial forces
    \begin{subequations}
    \label{FORCE:INERTIAL}
    \begin{alignat}{2}
        \label{FORCE:INERTIAL:B}
        F^\mathrm{bulk}_\mathrm{inertial}
        &= - \big( \rho\,\delt \uu + \Div(\uu\otimes\rho\uu) \uu \big)
        &&\tin\QT,
        \\
        \label{FORCE:INERTIAL:ST}
        F^{\mathrm{surf},\tau}_\mathrm{inertial}
        &= 0
        &&\ton\ST,
        \\
        \label{FORCE:INERTIAL:SN}
        F^{\mathrm{surf},\n}_\mathrm{inertial}
        &= 0
        &&\ton\ST,
    \end{alignat}
    \end{subequations}
    which are the forces resulting from the first variation of $\cA^\mathrm{kin}$,
    and the conservative forces
    \begin{subequations}
    \label{FORCE:CONS}
    \begin{alignat}{2}
        \label{FORCE:CONS:B}
        F^\mathrm{bulk}_\mathrm{cons}
        &= - \eps\, \Div(\Grad\phi \otimes \Grad\phi) - \Grad p
        &&\tin\QT,
        \\
        \label{FORCE:CONS:ST}
        F^{\mathrm{surf},\tau}_\mathrm{cons}
        &= \big[\eps\, (\Grad\phi \otimes \Grad\phi)\n - \delta\, \Divg(\Gradg\psi \otimes \Gradg\psi)\big]_\tau - \Gradg q
        &&\ton\ST,
        \\
        \label{FORCE:CONS:SN}
        F^{\mathrm{surf},\n}_\mathrm{cons}
        &= \big[ \eps\, (\Grad\phi \otimes \Grad\phi)\n \cdot \n \big] \n + p\n - qH\n
        &&\ton\ST,
    \end{alignat}
    \end{subequations}
    which result from the first variations of $\cA^\mathrm{bulk}$ and $\cA^\mathrm{surf}$.
    Here, $p:\QT\to\R$ and $q:\ST\to\R$ are Lagrange multipliers that account for the incompressibility constraints \eqref{COND:Y:B} and \eqref{COND:Y:S}. They can be interpreted as the bulk and the surface pressure, respectively. Notice that, because of \eqref{COND:Y:B} and \eqref{COND:Y:S}, we get
    \begin{align*}
        0 &= \intO p\, \Div(\tilde \y) \dx
        = - \intO \Grad p \cdot \tilde \y \dx 
            + \intG p\n \cdot \tilde \y \dHH{d-1},
        \\
        0 &= \intG q\, \Divg(\tilde \y) \dHH{d-1}
        = - \intG \Gradg q \cdot \tilde \y \dHH{d-1} 
            - \intG qH\n \cdot \tilde \y \dHH{d-1}
    \end{align*}
    via integration by parts.
    This explains why $p$ and $q$ appear in the equations \eqref{FORCE:CONS} in this particular way.

    \medskip

    \begin{remark}
        In the situation of unmatched densities ($\trho_1\neq \trho_2$), we expect to obtain the inertial forces
        \begin{alignat*}{2}
            F^\mathrm{bulk}_\mathrm{inertial}
            &= - \Big[ \delt\big( \rho(\phi)\,\uu\big) 
                + \Div\big(\uu \otimes \big(\rho(\phi)\uu + \J \big)\big) \uu 
            \Big]
            &&\tin\QT,
            \\
            F^{\mathrm{surf},\tau}_\mathrm{inertial}
            &= \frac 12 (\J\cdot\n) \uu_\tau
            &&\ton\ST,
            \\
            F^{\mathrm{surf},\n}_\mathrm{inertial}
            &= \frac 12 (\J\cdot\n) (\uu\cdot\n) \n
            &&\ton\ST.
        \end{alignat*}
        However, it remains unknown whether these identities can be derived from the least action principle. Following \cite{Liu2015}, it may still be possible to derive the Navier–Stokes equation in its weak form. However, this formulation would not allow for a proper identification of the inertial forces. 
    \end{remark}

    \subsubsection{Derivation of the dissipative forces via Onsager's principle of maximal energy dissipation}

    According to \textit{Onsager's principle of maximal energy dissipation}, the dissipative forces can be obtained via the first variation of the \textit{Rayleigh dissipation functional} 
    \begin{align*}
        \mathcal{R} \coloneqq \frac 12 \mathcal D
    \end{align*}
    with respect to the velocity field $\uu$, where $\mathcal D$ is given by \eqref{DEF:DISS:GLOB}.
    Let $\vv:\QT\to\R^d$ be a smooth vector field, which satisfies
    \begin{alignat}{2}
        \label{COND:V:B}
        \Div(\vv) &= 0 
        &&\ton\QT,
        \\
        \label{COND:V:S}
        \Divg(\vv) &= 0
        &&\ton\ST.
    \end{alignat}
    Then a straightforward computation yields
    \begin{align*}
        &\partial_\uu \mathcal{R}(\uu)[\vv]
        = \frac{\mathrm{d}}{\mathrm{d} \varsigma}\Big|_{\varsigma=0} \left[ \frac 12 \mathcal{D}(\uu +\varsigma\vv) \right]
        \\
        &\quad = \intO 2 \nu(\phi) \D \uu : \D\vv \dx
			+ \intG \gamma_{\tau} \uu_\tau \cdot \vv_\tau \dHH{d-1} 
			+ \intG \gamma_{\n} (\uu \cdot \n)\, \n \cdot \vv \dHH{d-1}
        \\
        &\quad = - \intO \Div\big( 2 \nu(\phi) \D \uu \big) \cdot \vv \dx
			+ \intG \big[ 2 \nu(\phi) (\D \uu) \n + \gamma_\tau \uu \big]_\tau \cdot \vv_\tau \dHH{d-1} 
        \\
        &\quad\qquad
			+ \intG \gamma_{\n} (\uu \cdot \n)\, \n \cdot \vv \dHH{d-1}.
    \end{align*}
    Consequently, the dissipative forces can be identified as
    \begin{subequations}
    \begin{alignat}{2}
        \label{FORCE:DISS:B}
        F^\mathrm{bulk}_\mathrm{diss}
        &= - \Div\big( 2 \nu(\phi) \D \uu \big)
        &&\tin\QT,
        \\
        \label{FORCE:DISS:ST}
        F^{\mathrm{surf},\tau}_\mathrm{diss}
        &= \big[ 2 \nu(\phi) (\D \uu) \n + \gamma_\tau \uu \big]_\tau
        &&\ton\ST,
        \\
        \label{FORCE:DISS:SN}
        F^{\mathrm{surf},\n}_\mathrm{diss}
        &= \gamma_{\n} (\uu \cdot \n)\, \n
        &&\ton\ST.
    \end{alignat}
    \end{subequations}

    \subsubsection{Force balances}

    According to Newton's second law of motion, the sum of inertial forces, conservative forces and dissipative forces is zero.
    More precisely, this means that 
    \begin{subequations}
    \begin{alignat}{5}
        \label{FB:B}
        &F^\mathrm{bulk}_\mathrm{inertial} 
        && + F^\mathrm{bulk}_\mathrm{cons} 
        && + F^\mathrm{bulk}_\mathrm{diss}
        &&= 0 
        &&\tin\QT,
        \\
        \label{FB:ST}
        &F^{\mathrm{surf},\tau}_\mathrm{inertial} 
        && + F^{\mathrm{surf},\tau}_\mathrm{cons} 
        && + F^{\mathrm{surf},\tau}_\mathrm{diss}
        &&= 0 
        &&\ton\ST,
        \\
        \label{FB:SN}
        &F^{\mathrm{surf},\n}_\mathrm{inertial} 
        && + F^{\mathrm{surf},\n}_\mathrm{cons} 
        && + F^{\mathrm{surf},\n}_\mathrm{diss}
        &&= 0 
        &&\ton\ST.
    \end{alignat}
    \end{subequations}

    In view of \eqref{FORCE:INERTIAL:B}, \eqref{FORCE:CONS:B} and \eqref{FORCE:DISS:B}, the force balance \eqref{FB:B} can be reformulated as
    \begin{align}
        \label{PDE:NS*}
        \delt (\rho \uu) + \Div(\uu\otimes\rho\uu) 
        = \Div\big( 
                2\nu(\phi)\, \D\uu
                - p \mathbf{I}
                - \eps \Grad\phi \otimes \Grad\phi
            \big)
        \quad\tin\QT.
    \end{align}
    Recalling \eqref{BL:CLM}, \eqref{DEF:TT}, and \eqref{SIMP:JRHO}, the stress tensor $\T$ can be identified as
    \begin{align}
		\label{DEF:T:E}
        \T 
        = 2 \nu(\phi) \D\uu 
        - p \I 
        - \eps \Grad\phi \otimes \Grad\phi
        \quad\tin\QT,
	\end{align}
    and equation \eqref{PDE:NS*} can be reformulated as
    \begin{align}
        \label{PDE:NS:E}
        \delt (\rho \uu) + \Div(\uu\otimes\rho\uu) 
        = \Div(\T)
        \quad\tin\QT.
    \end{align}

    Next, collecting \eqref{FORCE:INERTIAL:ST}, \eqref{FORCE:CONS:ST}, and \eqref{FORCE:DISS:ST}, and invoking \eqref{DEF:T:E}, the force balance \eqref{FB:ST} can be expressed as
    \begin{align}
        \label{BC:NAV:E}
        \big[ 
            \T \n 
            + \Gradg q
            + \gamma_\tau \uu 
        \big]_\tau
        =
        \big[ 
            - \delta\, \Divg(\Gradg\psi \otimes \Gradg\psi) 
        \big]_\tau
        \quad\ton\ST.
    \end{align}

    Lastly, recalling \eqref{DEF:T:E} and combining \eqref{FORCE:INERTIAL:SN}, \eqref{FORCE:CONS:SN}, and \eqref{FORCE:DISS:SN}, the force balance \eqref{FB:SN} can be reformulated as
    \begin{align*}
        \big[\T\n\cdot\n + \gamma_\n  (\uu\cdot\n) \big] \n 
			= - qH \n
        \ton\ST.
    \end{align*}
    In view of \eqref{BC:NORVEL}, this means that
    \begin{align}
        \label{BC:NOR:E}
        \T\n\cdot\n + \gamma_\n \cV 
			= - qH 
        \ton\ST.
    \end{align}

    \subsubsection{Determination of flux terms and boundary conditions}

    The last step is to complete the model derivation by identifying the flux terms $\J_\phi$ and $\J_\psi$ as well as the missing boundary conditions via the global energy dissipation law \eqref{EN:DISS:GLOB}. 
    To this end, we first need to compute the time derivative of the energy functional.
    Applying the Reynolds transport theorem for evolving domains (see, e.g., \cite[Theorem~2.11.1]{BDGP-book}) and integrating by parts, we calculate
    \begin{align*}
        &\ddt \cE^\mathrm{kin}(\uu)
        = \intO \deltb \left( \frac \rho 2 \abs{\uu}^2 \right) \dx
        = \intO \rho\, \deltb \uu \cdot \uu \dx
        \\[1ex]
        &\quad = \intO \big( \rho\,\delt \uu + \Div(\uu\otimes\rho\uu) \big)  
            \cdot \uu \dx
        = \intO \Div(\T)  
            \cdot \uu \dx
        \\[1ex]
        &\quad = 
        \intG \big[ \T\n \big]_\tau \cdot \uu_\tau \dHH{d-1}
        + \intG \big[ (\T\n)\cdot\n \big] (\uu\cdot\n) \dHH{d-1}
        - \intO \T : \Grad\uu \dx\,.
    \end{align*}
    Hence, using \eqref{DEF:T:E}--\eqref{BC:NOR:E}, this time derivative can be reformulated as
    \begin{align*}
        \ddt \cE^\mathrm{kin}(\uu)
        &= - \intG \big[ \Gradg q + \gamma_\tau \uu + \delta\, 
        \Divg(\Gradg\psi\otimes\Gradg\psi)\big]_\tau \cdot \uu_\tau \dHH{d-1}
        \\
        &\qquad 
        - \intG \big[ \gamma_\n (\uu\cdot\n) + qH \big] (\uu\cdot\n) \dHH{d-1}
        \\
        &\qquad
        - \intO \big[ 
            2 \nu(\phi) \D\uu 
            - p \I 
            - \eps \Grad\phi \otimes \Grad\phi 
        \big] : \Grad\uu \dx\,.
    \end{align*}
    Consequently, integrating by parts and recalling the identities $\Div \uu = 0$ in $\QT$ (cf.~\eqref{PDE:DIVV}) and $\Divg \uu_\tau - H (\uu \cdot \n) = \Divg \uu = 0$ on $\ST$ (cf.~\eqref{PDE:DIVV:S}), we conclude that
    \begin{align}
        \label{TDER:EKIN}
        \begin{split}
        \ddt \cE^\mathrm{kin}(\uu)
        &= 
        - \intO 2 \nu(\phi) \abs{\D \uu}^2 \dx
        - \intG \gamma_\tau \abs{\uu_\tau}^2 \dHH{d-1}
        - \intG \gamma_\n \abs{\uu_\n}^2 \dHH{d-1} 
        \\
        &\qquad 
        + \intO \eps (\Grad\phi \otimes \Grad\phi) : \Grad\uu \dx
        + \intG \delta (\Gradg\psi\otimes\Gradg\psi) : \Gradg\uu_\tau \dHH{d-1}\,.
        \end{split}
    \end{align}

    To compute the time derivatives of the free energy contributions, we use the identities
    \begin{align} 
        \deltb \Grad \phi - \Grad \deltb \phi  
        &= - (\Grad\uu)^T \Grad\phi
        \\
        \deltc \Gradg \psi - \Gradg \deltc \psi   
        &= \Gradg \uu  \Gradg \psi - \big[ \Gradg \uu + (\Gradg \uu)^T \big]_\tau \Gradg \psi,
    \end{align} 
    see, e.g., \cite[Lemma~38]{BGN-book}. 
    Using the Reynolds transport theorem for evolving domains (see, e.g., \cite[Theorem~2.11.1]{BDGP-book}), the transport theorem for evolving hypersurfaces (see, e.g., \cite[Theorem~2.10.1]{BDGP-book}), and integration by parts, we deduce that
    \begin{align}
        \label{TDER:FREE}
		&\ddt \big[ \cE^\mathrm{bulk}(\phi) + \cE^\mathrm{surf}(\phi,\psi) \big]
		= \intO \deltb f \dx + \intG \deltc g \dHH{d-1}
		\nonumber \\[1ex]
		&\begin{aligned}
            &= \intO \eps \Grad\phi \cdot \deltb\Grad\phi + \frac{1}{\eps} F'(\phi)\deltb \phi\dx 
    		\\
    		&\quad + \intG \delta \Gradg\psi \cdot \deltc\Gradg\psi + \frac{1}{\delta} G'(\psi) \deltc \psi  
    			+ \chi(K) (\alpha\psi-\phi) (\alpha \deltc \psi - \deltb \phi)\dHH{d-1}
        \end{aligned}
		\nonumber \\[1ex]
        &\begin{aligned}
    		& = \intO \eps \Grad\phi \cdot \Grad\deltb\phi - \eps \Grad\phi \cdot \big[(\Grad\uu)^T \Grad\phi \big] 
    			+ \frac{1}{\eps} F'(\phi)\deltb \phi\dx 
    		\\
    		&\quad + \intG \delta \Gradg\psi \cdot \Gradg\deltc\psi 
    			+ \delta \Gradg\psi \cdot \Big( \Gradg \uu \, \Gradg\psi 
                - \big[\Gradg \uu + (\Gradg \uu)^T\big]_\tau \Gradg\psi \Big) \dHH{d-1} 
            \\
            &\quad + \intG \frac{1}{\delta} G'(\psi) \deltc \psi 
                + \chi(K) (\alpha\psi-\phi) (\alpha \deltc \psi - \deltb \phi)\dHH{d-1} 
		\end{aligned}
		\nonumber \\[1ex]
        &\begin{aligned}
    		& = \intO \Big( - \eps \Lap \phi + \frac{1}{\eps} F'(\phi) \Big) \deltb \phi
    		- \eps (\Grad\phi \otimes \Grad\phi) : \Grad\uu \dx 
    		\\
    		&\quad + \intG \Big( -\delta \Lapg\psi + \frac{1}{\delta} G'(\psi) + \alpha\eps \deln\phi \Big) \deltc \psi  
    		- \delta (\Gradg\psi\otimes\Gradg\psi) : \Gradg\uu_\tau \dHH{d-1}
    		\\
    		&\quad + \intG \Big( \chi(K) (\alpha\psi-\phi) - \eps\deln\phi \Big) (\alpha\deltc\psi - \deltb\phi)\dHH{d-1}\,.
        \end{aligned}
	\end{align}
    Now, we introduce the chemical potentials
    \begin{alignat}{2}
        \label{PDE:MU:E}
        \mu &\coloneqq - \eps \Lap \phi + \frac{1}{\eps} F'(\phi)
        &&\quad\tin\QT,
        \\
        \label{PDE:THETA:E}
        \theta &\coloneqq -\delta \Lapg\psi + \frac{1}{\delta} G'(\psi) + \alpha\eps \deln\phi
        &&\quad\ton\ST.
    \end{alignat}
    Furthermore, recalling \eqref{EQ:PHI} and \eqref{EQ:PSI}, we have
    \begin{alignat}{2}
        \label{EQ:PHI:E}
        \deltb\phi + \Div(\J_\phi) &= 0
        &&\quad\tin\QT,
        \\
        \label{EQ:PSI:E}
        \deltc\psi + \Divg(\K_\phi) &= 0
        &&\quad\ton\ST.
    \end{alignat}
    Hence, combining \eqref{TDER:EKIN} and \eqref{TDER:FREE}, and recalling \eqref{ASS:NOFLUX} when integrating by parts, the time derivative of the total energy can be expressed as
    \begin{align}
        \label{TDER:ETOT}
        \begin{split}
        \ddt \cE
        &= 
        - \intO 2 \nu(\phi) \abs{\D \uu}^2 \dx
        - \intG \gamma_\tau \abs{\uu_\tau}^2 \dHH{d-1}
        - \intG \gamma_\n \abs{\uu_\n}^2 \dHH{d-1} 
         \\
        &\qquad 
        - \intO \mu\, \Div(\J_\phi) \dx
        - \intG \theta\, \Divg(\K_\psi) \dHH{d-1}
        \\
    	&\qquad + \intG \Big( \chi(K) (\alpha\psi-\phi) - \eps\deln\phi \Big) (\alpha\deltc\psi - \deltb\phi)\dHH{d-1}
        \end{split}
        \nonumber \\
        \begin{split}
        &= 
        - \intO 2 \nu(\phi) \abs{\D \uu}^2 \dx
        - \intG \gamma_\tau \abs{\uu_\tau}^2 \dHH{d-1}
        - \intG \gamma_\n \abs{\uu_\n}^2 \dHH{d-1} 
         \\
        &\qquad 
        - \intO \Grad \mu\cdot \J_\phi \dx
        - \intG \Gradg \theta\cdot \K_\psi \dHH{d-1}
        \\
    	&\qquad + \intG \Big( \chi(K) (\alpha\psi-\phi) - \eps\deln\phi \Big) 
            (\alpha\deltc\psi - \deltb\phi)\dHH{d-1}
        \end{split}
        \nonumber \\
        \begin{split}
        &= 
        - \,\mathcal{D}\, 
        + \intO \frac{1}{m_{\Omega}(\phi)}\Big( \J_\phi - m_{\Omega}(\phi)\, \Grad \mu \Big) \cdot \J_\phi \dx
        \\
        &\qquad + \intG \frac{1}{m_{\Gamma}(\psi)}\Big( \K_\psi - m_{\Gamma}(\psi)\, \Gradg \theta \Big) \cdot \K_\psi \dHH{d-1}
        \\
    	&\qquad + \intG \Big( \chi(K) (\alpha\psi-\phi) - \eps\deln\phi \Big) 
            (\alpha\deltc\psi - \deltb\phi)\dHH{d-1} \,.
        \end{split}
    \end{align}
    
    This means that the global energy dissipation law \eqref{EN:DISS:GLOB} is fulfilled if the relations
    \begin{alignat}{2}
        \label{EQ:JPHI:E}
        &\,\,\J_\phi = - m_{\Omega}(\phi)\, \Grad \mu
        &&\quad\tin\QT,
        \\
        \label{EQ:KPSI:E}
        &\K_\psi = - m_{\Gamma}(\psi)\, \Gradg \theta
        &&\quad\ton\ST,
        \\
        \label{BC:PHIPSI:E}
        &\begin{cases}
            \phi = \alpha \psi  &\text{if $K=0$},\\
            K \varepsilon \deln \phi = \alpha \psi - \phi &\text{if $K\in (0,\infty)$},\\
            \deln \phi = 0 &\text{if $K=\infty$}
        \end{cases}
        &&\quad\ton\ST
    \end{alignat}
    hold. In particular, because of the no-flux assumption \eqref{ASS:NOFLUX}, equation \eqref{EQ:JPHI:E} readily implies that the boundary condition 
    \begin{align}
        \label{BC:MUTHETA:E}
        m_\Omega(\phi) \deln \mu = 0 \quad\ton\ST
    \end{align}
    must be fulfilled. This corresponds to the choice $L=\infty$ in system \eqref{eqs:NSCH}.
    Finally, plugging \eqref{EQ:JPHI:E} and \eqref{EQ:KPSI:E} into \eqref{EQ:PHI:E} and \eqref{EQ:PSI:E}, respectively, we conclude that
    \begin{alignat}{2}
        \label{PDE:PHI:E}
        \deltb\phi &= \Div\big( m_{\Omega}(\phi)\, \Grad \mu \big)
        &&\quad\tin\QT,
        \\
        \label{PDE:PSI:E}
        \deltc\psi &= \Divg\big( m_{\Gamma}(\psi)\, \Gradg \theta \big)
        &&\quad\ton\ST.
    \end{alignat}

    In summary, for $L=\infty$ and the situation of matched densities $\trho_1=\trho_2$ (which entails $\J=\mathbf{0}$ in $\QT$), we have derived system \eqref{eqs:NSCH} along with formulae \eqref{eq:rho} and \eqref{eq:T} with the following correspondences:
    \begin{center}
    \makeatletter \setlength{\tabcolsep}{2pt} \makeatother
    \begin{tabular}{rclcrclcrcl}
        \eqref{eqs:NSCH:1} 
        &$\widehat{=}$ 
        &\eqref{PDE:NS:E}, \eqref{PDE:DIVV}, 
        &$\phantom{xx}$ 
        &\eqref{eqs:NSCH:2} 
        &$\widehat{=}$ 
        &\eqref{PDE:PHI:E}, 
        &$\phantom{xx}$ 
        &\eqref{eqs:NSCH:3} 
        &$\widehat{=}$ 
        &\eqref{PDE:MU:E}, 
        \\
        \eqref{eqs:NSCH:4} 
        &$\widehat{=}$ 
        &\eqref{PDE:PSI:E},
        &$\phantom{xx}$ 
        &\eqref{eqs:NSCH:5} 
        &$\widehat{=}$ 
        &\eqref{PDE:THETA:E},
        &$\phantom{xx}$ 
        &\eqref{eqs:NSCH:6} 
        &$\widehat{=}$ 
        &\eqref{BC:PHIPSI:E}, \eqref{BC:MUTHETA:E}, 
        \\
        \eqref{eqs:NSCH:7} 
        &$\widehat{=}$ 
        &\eqref{BC:NAV:E}, 
        &$\phantom{xx}$ 
        &\eqref{eqs:NSCH:8} 
        &$\widehat{=}$ 
        &\eqref{BC:NOR:E},
        &$\phantom{xx}$ 
        &\eqref{eqs:NSCH:9} 
        &$\widehat{=}$ 
        &\eqref{PDE:DIVV:S}, \eqref{BC:NORVEL},
        \\
        \eqref{eq:rho} 
        &$\widehat{=}$ 
        &\eqref{DEF:RHO*},
        &$\phantom{xx}$ 
        &\eqref{eq:T} 
        &$\widehat{=}$ 
        &\eqref{DEF:T:E}, \eqref{DEF:S}.
    \end{tabular}
    \end{center}
    Therefore, our model derivation is complete.

    \section{Conclusion}
    
    In this work, we have derived a thermodynamically consistent Navier–Stokes–Cahn–Hilliard model for two-phase flows in an evolving domain with a free boundary. The model incorporates bulk–surface interactions, variable contact angles, and a generalized Navier slip boundary condition, thereby extending previous diffuse-interface formulations to situations involving a domain evolution driven by the fluid motion. A central feature of our framework is that both the bulk and surface phase-fields evolve according to a coupled convective bulk–surface Cahn–Hilliard system, which accounts for material transfer between the bulk and the boundary as well as a variable contact angle.
    
    Starting from local mass balance laws, we systematically derived the governing equations and boundary conditions in two different ways: by the Lagrange Multiplier Approach and, assuming matched densities and no mass flux between bulk and surface, by the Energetic Variational Approach. The resulting system satisfies the fundamental physical principles of mass conservation, local energy dissipation, and preservation of the bulk volume and surface area of the evolving domain. We further showed that several classical and recently proposed models (cf.\cite{Abels2012,Giorgini2023}) appear as special limiting cases of our formulation. This highlights the generality and flexibility of the present model.
    
    The free boundary setting considered here opens up new perspectives for modeling interfacial phenomena involving deformable substrates, biological membranes, or motile cells. Future work shall address analytical questions such as well-posedness as well as numerical methods capable of approximating the coupled bulk–surface dynamics on evolving geometries.
    Moreover, applications to concrete physical or biological systems (e.g., droplet spreading on soft surfaces or intracellular flows within moving cells) offer promising directions for further study.

    \pagebreak[4]

    \section*{Acknowledgement}
    We thank the anonymous referees for their careful reading and valuable suggestions, which helped us to improve the manuscript.
    Both authors were partially supported by the ``RTG 2339'' \textit{IntComSin} funded by the Deutsche Forschungsgemeinschaft (DFG, German Research Foundation) -- Project-ID 321821685. 
    PK was additionally supported by the Deutsche Forschungsgemeinschaft (DFG, German Research Foundation) – Project-ID 52469428. YL was partially supported by the startup funding from Nanjing Normal University, the Natural Science Foundation of Jiangsu Province (Grant No. BK20240572), the Natural Science Foundation of the Jiangsu Higher Education Institutions of China (Grant No. 24KJB110020), and the China Postdoctoral Science Foundation under Grant Number 2025M773078.
    The support is gratefully acknowledged.
    	
    \section*{Compliance with Ethical Standards}
    
    \noindent\textbf{Date avability.}
    Data sharing not applicable to this article as no datasets were generated during the current study.
    
    \noindent\textbf{Conflict of interest.}
    The authors declare that there are no conflicts of interest.


    \footnotesize
    \bibliographystyle{abbrv}
    \bibliography{KL_NSCH}
    \normalsize


\end{document}